\documentclass[12pt]{amsart}

\usepackage{amsthm}
\usepackage{graphicx}
\input{amssym.def}
\input{amssym}
\input{epsf}
\usepackage{cite}
\def\R{\mathbb R}
\def\N{\mathbb N}
\def\T{\mathbb T}

\def\eps{\varepsilon}

\def\eps{\varepsilon}

\def\tilde{\widetilde}

\def\2dr#1#2{\left. \frac{d^2}{d{#1}^2} \right |_{#2}}
\def\d2#1{\frac{d^2}{d{#1}^2}}

\def\sm{\smallskip}

\newtheorem{theorem}{Theorem}
\newtheorem{lemma}[theorem]{Lemma} %unique numerotation
\newtheorem{proposition}[theorem]{Proposition}

\theoremstyle{definition}
\newtheorem{definition}[theorem]{Definition}
\newtheorem*{remark}{Remark}

\newcommand{\di}{d}
\newcommand{\dii}{{d-1}}
\newcommand{\ir}{\int_{\R^\di}}
\newcommand{\is}{\int_{S^{\dii}}}
\newcommand{\ph}{\varphi}
\renewcommand{\th}{\vartheta}
\newcommand{\TD}{{[0,T]}}
\newcommand{\XD}{{\T^\di}}

\newcommand{\sign}{\text{sign}}
\DeclareMathOperator*{\limm}{\overline{\lim}}

\def\Ack{\medskip\noindent {\bf Acknowledgements:}\ \ignorespaces}

\def\signig{\bigskip \begin{center} 
{\sc Irene M. Gamba\par\vspace{3mm}
Department of Mathematics, University of Texas at Austin \\
Austin, TX 78712 U.S.A. 
\par\vspace{3mm}
e-mail:} \tt{gamba@math.utexas.edu}  \end{center}}

\def\signvp{\bigskip \begin{center} 
{\sc Vladislav Panferov\par\vspace{3mm}
Department of Mathematics and Statistics \\ McMaster University, 
1280 Main St. West \\Hamilton, ON L8S 4K1 Canada
\par\vspace{3mm}
e-mail:} \tt{panferov@math.mcmaster.ca}  \end{center}}

\def\signcv{\bigskip \begin{center} 
{\sc C\'edric Villani\par\vspace{3mm} 
UMPA, ENS Lyon,
46 all\'ee d'Italie\par 
69364 Lyon Cedex 07,
France\par\vspace{3mm}
e-mail:} \tt{cvillani@umpa.ens-lyon.fr} \end{center}}

\begin{document}
\title[Maxwellian bounds for the Boltzmann equation]
{Upper Maxwellian bounds for the spatially homogeneous 
Boltzmann equation}

\author{I. M. Gamba}
\author{V. Panferov}
\author{C. Villani}

\begin{abstract}
For the spatially homogeneous Boltzmann equation with cutoff 
hard potentials it is shown that solutions remain bounded 
from above, uniformly in time, by a Maxwellian distribution, 
provided the initial data have a Maxwellian upper bound. 
The main technique is based on a comparison principle that 
uses a certain dissipative property of the linear Boltzmann 
equation. Implications of the technique to propagation 
of upper Maxwellian bounds in the spatially-inhomogeneous 
case are discussed. 
%
%We consider solutions of the spatially homogeneous Boltzmann 
%equation for hard potentials with an angular cutoff. 
%We establish that if the initial datum is bounded from above 
%by a Maxwellian distribution, the solution remains bounded 
%for all times by another time-independent Maxwellian 
%distribution. Our main technique is based on a comparison 
%principle which uses a certain monotonicity 
%property of the linear Boltzmann equation. 
\end{abstract}

\maketitle
\vspace{-0.5cm}
\section{Introduction and main result}

The nonlinear Boltzmann equation is a classical model 
for a gas at low or moderate densities. The gas in a 
spatial domain \(\Omega\subseteq\R^\di\), \(\di\ge 2\),
is modeled by the mass density function \(f(x,v,t)\), 
\((x,v)\in \Omega\times\R^\di\), where \(v\) is the velocity 
variable, and \(t\in \R\) is time. The equation 
for \(f\) reads 
\begin{equation}
\label{eq:boltz}
(\partial_t + v\cdot \nabla_x) f = Q(f)\,,
\end{equation}
where  \(Q(f)\) is a quadratic integral operator,
expressing the change of \(f\) due to instantaneous 
binary collisions of particles. The precise form 
of \(Q(f)\) will be introduced below, cf. 
also~\cite{CeIlPu,Vi}.  

Although some of our results deal with more general situations,
we will be mostly concerned with a special class of 
solutions that are independent of the spatial variable
(spatially homogeneous solutions). In this case \(f=f(v,t)\)
and one can study the initial-value problem
\begin{equation}
\label{eq:bosh}
\partial_t f = Q(f),\quad f|_{t=0}=f_0, 
\end{equation}
where \(0\le f_0\in L^1(\R^\di)\).
The spatially homogeneous theory is very well developed although not
complete. In the present paper we shall solve one of the most
noticeable open problems remaining in the field, by establishing the
following result.

\begin{theorem}
\label{thm:1}
Assume that \(0\le f_0(v) \le M_0(v)\), for a.\;a. 
\(v\in \R^\di\), where \(M_0(v) = e^{-a_0|v|^2+c_0}\) 
is the density of a Maxwellian distribution, \(a_0>0\), 
\(c_0\in\R\). Let \(f(v,t)\), \(v\in\R^\di\), \(t\ge0\) 
be the unique solution of equation~\eqref{eq:bosh} for 
hard potentials with the angular 
cutoff assumptions~\eqref{eq:asmp_b},~\eqref{eq:asmp_h}, 
that preserves the initial mass and energy~\eqref{eq:me}. 
Then there are constants \(a>0\) and \(c\in\R\) 
such that  \(f(v,t)\le M(v)\), for a.\;a. 
\(v\in \R^\di\) and for all \(t\ge 0\), 
where \(M(v)= e^{-a|v|^2+c}\). 
\end{theorem} 

Before going on, let us make a few comments about the interest
of these bounds.
Maxwellian functions 
\[
M(v) = e^{-a|v|^2+b\cdot v + c},\quad\text{with}\;\;
a>0, \,c\in \R, \;\,b\in\R^\di\text{ constants},
\]
are unique, within integrable functions, equilibrium 
solutions of~\eqref{eq:bosh}, and they 
provide global attractors for the time-evolution 
described by~\eqref{eq:bosh} (or~\eqref{eq:boltz}, 
with appropriate boundary conditions). 
Classes of functions bounded above by Maxwellians 
provide a convenient analytical framework 
for the local existence theory of strong solutions 
for~\eqref{eq:boltz}, see Grad~\cite{Gr} 
and Kaniel-Shinbrot~\cite{KaSh}. Such bounds 
also play an important role in the proof of 
validation of the Boltzmann equation 
by Lanford~\cite{La}, see also~\cite{CeIlPu}. 
However, establishing the propagation of uniform 
bounds is generally a difficult problem, solved 
only in the context of small solutions in an 
unbounded space, see Illner-Shinbrot~\cite{IlSh} 
and subsequent works~\cite{BeTo,Ha,MiPe,Go}. 
The above results rely in a crucial way on the decay of 
solutions for large \(|x|\) and on the dispersive 
effect of the transport term, in order to control 
the nonlinearity. These 
effects may not be significant in other physical situations, 
and the spatially homogeneous problem presents 
a simplest example of such regime. 

In the spatially homogeneous case many additional 
properties of solutions can be established. Upper 
bounds with  polynomial decay for \(|v|\) 
large hold uniformly in time, see
Carleman~\cite{Ca0,Ca} and Arkeryd~\cite{Ar1}. 
Solutions are also known to 
have a lower Maxwellian bound for all positive 
times, even for compactly supported 
initial data~\cite{PuWe}. Many results have been 
established that concern the behavior of the 
moments with respect to the velocity variable, 
following the work by Povzner~\cite{Po}, 
see in particular~\cite{Ar,El,De,MiWe,Bo}. 
The Carleman-type estimates~\cite{Ca0,Ca,Ar1}
were crucial in the treatment of the weakly inhomogeneous 
problem given in~\cite{ArEsPu}. 
However, as also pointed out in ref.~\cite{ArEsPu}, 
Maxwellian bounds of the local existence theory~\cite{Gr,KaSh}
are not known to hold on longer time-intervals, 
and it remains an open problem to characterize 
the approach to the Maxwellian equilibrium in classes of 
functions with exponential decay. 
%The problem of qualitative behavior of solutions 
%of the spatially homogeneous Boltzmann 
%equation~\eqref{eq:bosh} 
%has been extensively studied, starting from the 
%pioneering work by Carleman~\cite{Ca0,Ca}; for later 
%developments we refer to the works
%\cite{Ca,Po,Ar,Di,El,De,MiWe} and the reviews in 
%\cite{CeIlPu,Vi}.   
%The upper bounds stated in Theorem~\ref{thm:1}
%are very natural, since the Maxwellian functions 
%\[
%M(v) = e^{-a|v|^2+b\cdot v + c}\quad\text{with}\quad
%a>0, \;b\in\R^\di,\; c\in \R,
%\]
%are the only steady density solutions for the 
%problem~\eqref{eq:bosh}. Solutions in classes 
%of functions bounded above by Maxwellians also 
%appear in the contexts of short-time existence 
%for the spatially inhomogeneous problem~\cite{Gr,KaSh}, 
%rigorous validity of the Boltzmann 
%equation~\cite{La,CeIlPu}, global existence 
%for small data~\cite{IlSh,BeTo,Go}, perturbations 
%of the spatially-homogeneous states~\cite{ArEsPu}.
%The latter reference in particular points out 
%a gap in the theory of the 
%spatially-homogeneous problem, namely that the approach
%to equilibrium is not easily characterized in 
%classes of functions with Maxwellian tails. 
The present work aims to at least partially 
remedy this situation, and to develop a technique 
that could be used to obtain further results 
in this direction.  
%
%
%However, known results
%about the long-time behavior of solutions concern 
%the ``\(L^1\) behavior'', and the only uniform 
%in time pointwise estimates are those of polynomial 
%decay in the velocity space, due to 
%Carleman~\cite{Ca}.  
%
%The issue of propagation of uniform Maxwellian 
%estimates becomes even more important (and 
%significantly more difficult) in the spatially 
%inhomogeneous case, where it appears in the 
%results by Grad, Kaniel and Shinbrot, Lanford
%and Arkeryd, Esposito, Pulvirenti. While our 
%results do not apply directly to that case we believe 
%that some of the methods may be useful for 
%studying problems with spatial dependence.  
%

We will next introduce the notation and the necessary 
concepts to make the statement of Theorem~\ref{thm:1}
more precise. 
%
%We will return to the discussion of solutions 
%of~\eqref{eq:bosh}, but first we need to introduce 
%the particular form of the collision term \(Q(f)\)
%and to give a more precise meaning to the concepts 
%appearing in the formulation of Theorem~\ref{thm:1}. 
%
We set in~\eqref{eq:bosh}
%formulate the precise assumptions that appear 
%in the conditions of Theorem~\ref{thm:1}. We set
%
%The collision operator \(Q(f)\) appearing 
%in~\eqref{eq:boltz} and~\eqref{eq:bosh} is 
%as follows, in the case when \(f=f(v)\), 
\begin{equation}
\label{eq:co}
Q(f)\,(v,t)= \int_{\R^\di} \int_{S^{\di-1}} 
 (f_*'\,f'-f_*\,f )\,B(v-v_*,\sigma)\, d\sigma\,dv_*, 
\end{equation}
where, adopting common shorthand notations, 
\(f=f(v,t)\), \(f'=f(v',t)\), \(f_*=f(v_*,t)\), 
\(f'_*=f(v'_*,t)\). Here \(v\), \(v_*\) denote 
the velocities of two particles either before or 
after a collision, 
\begin{equation}
\label{eq:cp}
v' \,=\, \frac{v+v_*}2 
+ \frac{|v-v_*|}2\,\sigma\,, \qquad
v'_* \,=\, \frac{v+v_*}2 
- \frac{|v-v_*|}2\,\sigma\,,
\end{equation}
are the transformed velocities, 
and \(\sigma\in S^{\dii}\) is 
a parameter determining the direction of the relative 
velocity \(v'-v'_*\). 
In the more general case of~\eqref{eq:boltz}, the space
variable \(x\) appears (similarly to \(t\) above) 
in each occurrence of \(f\), \(f_*\), 
\(f'\), \(f'_*\); we shall often omit 
the \(t\) and \(x\) variables from the notation for brevity. 

Many properties of the solutions of the Boltzmann 
equation depend crucially on certain features of the
kernel \(B\) in \eqref{eq:co}. Its physical meaning 
is the product of the magnitude of the relative 
velocity by the effective scattering 
cross-section (see~\cite[\S 18]{LaLi1} for terminology 
and explicit examples); 
this quantity characterizes the relative frequency 
of collisions between particles. 
Our assumptions on \(B\) fall 
in the category of ``hard potentials with angular
cutoff'', cf.~\cite{Vi}.
More precisely, we assume that 
%\(B\) depends on \(|v-v_*|\) and 
%\(\cos\th = \tfrac{(v-v_*)\cdot\sigma}{|v-v_*|}\) 
%only, and it satisfies 
\begin{equation}
\label{eq:asmp_b}
B(v-v_*,\sigma) = |v-v_*|^\beta\, h(\cos\th), 
\quad \cos\th = \tfrac{(v-v_*)\cdot\sigma}{|v-v_*|},
\end{equation}
where \(0<\beta\le 1\) is a constant and 
\(h\) is a nonnegative function on \((-1,1)\) such that 
\begin{equation}
\label{eq:asmp_hi}
h(z)+h(-z) \quad\text{is nondecreasing on \((0,1)\)}
\end{equation}
and
\begin{equation}
\label{eq:asmp_h}
0\le h(\cos\th)\sin^\alpha\th \le C, \quad\th\in (0,\pi),
\end{equation}
where \(\alpha<\di-1\) and \(C\) is a constant. 
Assumption~\eqref{eq:asmp_h} implies in particular 
that the integral \(\is h(\cos\th)\, d\sigma\) is 
finite; for convenience we normalize it by setting
\begin{equation}
\label{eq:norm_h}
\is h(\cos\th)\, d\sigma = \omega_{\di-2} 
\int_{-1}^1 h(z)\,(1-z^2)^{\frac{\di-3}{2}} \, dz = 1,  
\end{equation}
where \(\omega_{\di-2}\) is the measure of the 
\((\di-2)\)-dimensional sphere. 
The classical hard-sphere model in \(\R^\di\), 
satisfies~\eqref{eq:asmp_b}
with \(\beta=1\), ~\eqref{eq:asmp_hi} 
and~\eqref{eq:asmp_h} with \(\alpha=\di-3\). 

%The growth of the kernel \(B\) as \(|v-v_*|\to\infty\)
%is essential for our methods. On the other hand, 
%it is conceivable that the condition on the angular 
%dependence~\eqref{eq:asmp_h} could be further 
%relaxed to include certain non-cutoff kernels. 
%We will make no attempt to do so here, however. 
%\begin{equation}
%\label{eq:asmp_B}
%\frac{C_2\,|v-v_*|\,(1+|v-v_*|)^{\beta-1}}{ (\sin\frac{\th}{2})^\alpha}\le 
%B(v-v_*,\sigma) \le \frac{C_1\,(1+|v-v_*|^\beta)}
%{(\sin\frac{\th}{2})^\alpha},
%\end{equation}
%for certain constants \(C_1\), \(C_2>0\) and 
%\(0<\beta\le 1\) and \(0\le\alpha< {\di-1}\). 
%In particular, the classical ``hard sphere'' 
%model satisfies these assumptions with \(\beta=1\)
%and \(\alpha=\di-3\). On the other hand, the 
%pseudo-Maxwell model is excluded by the condition
%\(\beta>0\).  
%
%The growth of the kernel
%for \(|v-v_*|\to\infty\) seems to be essential for 
%our result (notice that the condition \(\beta>0\)
%excludes the pseudo-Maxwell model). On the other 
%hand, it is possible that the condition 
%\(\alpha< {\di-1}\), implying the angular 
%cutoff, can be relaxed, allowing for certain 
%non-cutoff kernels. We have not tried to do 
%such an extension here. 

Notice that we can write 
\(
Q(f)= Q^+(f)-Q^-(f), 
\)
where \(Q^+(f)\) is the ``gain'' term, 
and \(Q^-(f)\) is the ``loss'' term, 
\[
Q^+(f) = \ir \is f' f'_*\,B(v-v_*,\sigma)\, 
d\sigma\,dv_*, \quad
Q^-(f) =  (f*|v|^\beta)\,f , 
\] 
and \(*\) denotes the convolution in \(v\). 
Because of the symmetry 
\(\sigma\mapsto -\sigma\) in the integral 
defining \(Q^+(f)\) we can restrict the 
\(\sigma\)-integration above  to the 
half-sphere \(\{\cos\th>0\}\) if we  
simultaneously replace \(B(v-v_*,\sigma)\) by 
\[
\overline{\!B}(v-v_*,\sigma) 
: = (B(v-v_*,\sigma)+B(v-v_*,-\sigma))\, 1_{\{\cos\th>0\}}.
\] 
It will be convenient to introduce the following 
(nonsymmetric) bilinear forms of the collision
terms,
\begin{equation}
\label{eq:bforms}
Q^+(f,g) = \ir \is f'_* \,g'\,\, 
\overline{\!B}(v-v_*,\sigma)\,
d\sigma\,dv_*,  \quad
Q^-(f,g) =  (f*|v|^\beta)\,g ,  
\end{equation}
for which obviously \(Q^\pm(f)=Q^\pm(f,f)\). 

We say that a nonnegative function \(f\in C([0,\infty); 
L^1(\R^\di))\), such that \((1+|v|^2)f\in L^\infty((0,\infty); 
L^1(\R^\di))\),
is a (mild) solution of~\eqref{eq:bosh} if  for 
almost all \(v\in \R^\di\)
\begin{equation}
\label{eq:mild}
f(v,0)=f_0(v); 
\qquad 
f(v,t)-f(v,s) = \int_s^t Q(f)(v,\tau)\, d\tau, 
\end{equation}
for all \(0\le s<t\). Notice that the conditions on \(f\)
imply (in the spatially-homogeneous case!) that 
\begin{equation}
\label{eq:ql1}
Q^+(f),\,Q^-(f)\in L^\infty((0,\infty); 
L^1(\R^\di)),
\end{equation}
so the integral form in~\eqref{eq:mild} 
is well-defined. This also implies that \(f\) is weakly 
differentiable with respect to \(t\) and that the 
differential equation~\eqref{eq:bosh} holds in the sense of 
distributions on \(\R^\di\times(0,\infty)\). 
 
The existence of a unique solution satisfying the 
conservations of mass and energy, 
\begin{equation}
\label{eq:me}
\ir f(v,t)\,dv = \ir f_0(v)\,dv, \qquad
\ir f(v,t)\, |v|^2\,dv = \ir f_0(v)\, |v|^2\,dv
\end{equation}
follows from a theorem by Mischler and 
Wennberg~\cite{MiWe}, for all \(f_0\ge 0\) for which the above
integrals are finite. The second condition in~\eqref{eq:me}
is also necessary for the uniqueness~\cite{We}.
For the initial data with strong decay (as in 
Theorem~\ref{thm:1}) one could also 
refer to the well-known results by Carleman, Arkeryd 
and DiBlasio~\cite{Ar,Ar1,Di}. 
The following theorem summarizes the main results 
about qualitative properties of solutions in the 
case of ``hard potentials with cutoff'' 
known before this work.

{\parindent=0pt
\parskip=10pt
\begin{theorem}
%[{\em Estimates for the spatially homogeneous Boltzmann equation with
%cut-off hard potentials}] 
\label{alreadyknown}
Let $f(v,t)$, \(v\in\R^\di\), \(t\ge 0\), (\(n\ge 2\))
be a solution of~\eqref{eq:bosh} that satisfies~\eqref{eq:me},
and let the kernel \(B\) in the Boltzmann operator~\eqref{eq:co} 
satisfy~\eqref{eq:asmp_b},~\eqref{eq:asmp_h}. Then
%\begin{enumerate}
%\addtolength{\itemsep}{0.5\baselineskip}
%\addtolength{\labelsep}{2pt}
%\addtolength{\leftmargin}{10pt}
%\addtolength{\rightmargin}{20pt}
%\addtolength{\itemindent}{-10pt}
%\setlength{\labelwidth}{1.5cm}
%\renewcommand{\labelenumi}{(\roman{enumi})}
%\item 

{\hspace{6pt}}\parbox{14.7cm}{
\parindent=0pt
\parskip=8pt
{\em (i)}\hspace{4pt}
if \(f_0\in L^\infty(\R^\di)\) then 
\(f(t,\cdot)\in L^\infty(\R^\di)\), \(t\ge0\). Moreover,  
if\; $(1+|v|)^{s}f_0\,\in L^\infty(\R^\di_v)$ for some $s>s_0$, then 
$(1+|v|)^{s}f(v,t)\in L^\infty(\R^\di_v)$, \(t\ge0\). Here \(s_0\)
is a constant dependent on the dimension \(\di\). 
}

%\item
{\hspace{6pt}}\parbox{14.7cm}{
\parindent=0pt
\parskip=8pt
{\em (ii)}\hspace{4pt}
if the integral of \(f\) is nonzero, then for every
$t_0>0$ there is a Maxwellian \(M(v)={K} e^{-\kappa|v|^2}\),
$K>0$, $\kappa>0$ such that
\[ 
f(v,t) \geq M(v),
\quad t\geq t_0,\quad 
\text{for a.\:a. }  v\in \R^\di. 
\]
}

%\item
{\hspace{6pt}}\parbox{14.7cm}{
\parindent=0pt
\parskip=8pt
{\em (iii)}\hspace{4pt}
for all $t_0>0$ and for all $k>1$, the quantity
\( m_k(t)= \ir f(v,t)\,|v|^{2k}\,dv\)
is bounded uniformly for $t\geq t_0$; moreover, this bound is uniform
in $t\geq 0$ if $m_k(0)<+\infty$.
}

%\item
{\hspace{6pt}}\parbox{14.7cm}{
\parindent=0pt
\parskip=8pt
{\em (iv)}\hspace{4pt}
In the case \(\di=3\) and \(B(v-v_*,\sigma)=c\,|v-v_*|\) 
(hard spheres) or \(B(v-v_*,\sigma)
=h(\frac{(v-v_*)\cdot\sigma}{|v-v_*|})\), 
\(h\in L^1(-1,1)\) (pseudo-Maxwell particles) if $f_0$ satisfies
\[ 
\frac{f_0}{M_0}\in L^1(\R^\di)
\]
for some Maxwellian \(M_0(v)=e^{-a_0|v|^2}\), $a_0>0$, 
then there exists constants $a>0$, \(C\) such that
\[ 
\ir \frac{f(v,t)}{M(v)}\,dv \le C,
\]
where \(M(v)=e^{-a|v|^2}\). 
}
%\end{enumerate}
\end{theorem}
}

Part (i) of this theorem is due to Carleman~\cite{Ca} 
in the case of the hard spheres; the general case was studied
by Arkeryd in~\cite{Ar1}. 
Part (ii) is due to A.~Pulvirenti and Wennberg~\cite{PuWe}. 
Part (iii) is due to Desvillettes~\cite{De} under the 
additional assumption that a moment \(m_{k_0}(t)\) 
of order \(k_0>1\) is finite initially; this assumption 
was removed by 
Mischler and Wennberg~\cite{MiWe}. Earlier result by 
Arkeryd~\cite{Ar} and Elmroth~\cite{El} state that all 
moments remain 
bounded uniformly in time, once they are finite initially.  
Finally, part (iv) is due to Bobylev~\cite{Bo}; 
we will give an extension of this result to the 
class of Boltzmann kernels 
satisfying~\eqref{eq:asmp_b}--\eqref{eq:asmp_h}
in Section~\ref{sec:mom}.  
%
%The proofs of all statements included in 
%Theorem~\ref{alreadyknown} are based on constructive 
%methods; all the constants appearing in the estimates 
%can be made explicit. Results (ii) and
%(iii) express properties which hold true essentially 
%independently of the
%initial datum, while (i) and (iv) express properties 
%of propagation of decay,
%in some sense. As a consequence of (i) and (iv), if the 
%initial data decays fast enough,
%then the solution is bounded and has finite square-exponential 
%moments.

Our main contribution in the present work is to show
that the estimates for the spatially homogeneous 
Boltzmann equation (precisely, parts (i) and (iv) 
of Theorem~\ref{alreadyknown}, together with the 
conservation of mass) imply Theorem~\ref{thm:1}. 
Since we do not use other properties of the 
spatially-homogeneous problem we can state our 
result in a more general, spatially inhomogeneous setting. 

We consider solutions 
of~\eqref{eq:boltz} with the spatial domain 
\(\Omega=\T^\di\) (\(\di\)-dimensional torus,
or the unit hypercube with periodic 
boundary conditions), 
on an arbitrary finite time interval \([0,T]\). 
Spatially homogeneous solutions 
are then a special subclass characterized by the 
constant dependence on the 
\(x\) variable. 
%We will impose apriori rather mild 
%restrictions on the spatial and temporal regularity 
%of the solutions. For the applications to 
%Theorem~\ref{thm:1} it is sufficient to require 
%that 
%\begin{equation}
%\label{eq:qint}
%f(x,v,t)\, f(x,v_*,t)\, B(v-v_*,\sigma)\in 
%L^1(\T^\di_x\times\R^\di_v\times\R^\di_{v_*}
%\times(0,T)\times S^\dii_\sigma),  
%\end{equation}
%and that~\eqref{eq:boltz} is satisfied in the 
%sense of distributions. 
%For the spatially homogeneous solutions,~\eqref{eq:qint}
%is an immediate consequence of our assumptions 
%on \(f\) and the kernel \(B\). In the 
%spatially inhomogeneous case condition~\eqref{eq:qint} can 
%be significantly relaxed, as indicated below.  
%
%Moreover, condition~\eqref{eq:qint} can 
%be replaced by requiring \(f\) to be a dissipative 
%solution of~\eqref{eq:boltz}, in the sense of 
%P.-L. Lions~\cite{Li} (see Definition~\ref{def:dis} 
%in Section~\ref{sec:cp}). 
%
To simplify the presentation, let us assume sufficient 
regularity (smoothness) of the solutions \(f(x,v,t)\) 
with respect to the \(x\) and \(t\) variables; 
this is not a restriction in the setting of Theorem 1, 
and the requirements of smoothness will be relaxed 
significantly later on to include a sufficiently 
wide class of weak solutions of the spatially 
inhomogeneous problem.
%
%Our result stated in Theorem~\ref{thm:1} depends in a crucial 
%way on the parts (i) and (iv) of Theorem~\ref{alreadyknown}
%(more precisely, the extended vesion of (iv) given in 
%Theorem~\ref{thm:int} below). 
%In fact these two properties are the main reason why our 
%result is restricted to the spatially homogeneous case. 
%Our main contribution in the present work is the following 
%conditional result which together with the parts (i) and 
%(iv) of Theorem~\ref{alreadyknown}  implies  Theorem~\ref{thm:1}. 

\begin{theorem}
%[{\em Conditional propagation of the Maxwellian upper bound}] 
\label{conditional}
Let $T>0$ and let $f\in C([0,T];L^1(\T^\di\times\R^\di))$, \(f\ge0\), 
be a (sufficiently regular) solution of the Boltzmann 
equation~\eqref{eq:boltz}, 
%satisfying~\eqref{eq:qint}, or, more generally, the dissipative 
%property of Definition~\ref{def:dis}, 
with the initial condition
\[
f(x,v,0)=f_0(x,v)\le M_0(v), \quad \text{for a. a. } 
(x,v)\in\T^\di\times\R^\di,
\] 
where \( M_0(v) = e^{-a_0 |v|^2 +c_0}\), \(a_0>0\), 
\(c_0\in\R\). Assume that the solution \(f(x,v,t)\) 
satisfies the estimates 
%\begeq\label{bound-f}
\begin{equation}
\label{eq:rho0}
\ir f(x,v,t)\, dv \,\ge\, \rho_0,\quad  
(x,t)\in \T^\di\times[0,T],
\end{equation}
and
\begin{equation}
\label{eq:fbds}
\sup_{(x,t)\,\in\,\T^\di\times[0,T]} \|f(x,v,t)\|_{L^\infty_v} \,\le\, C_0,
\qquad
\sup_{(x,t)\in \T^\di\times[0,T]} \ir \frac{f(x,v,t)}{M_1(v)}\,dv \leq
C_1, 
\end{equation}
%\endeq 
where \( M_1(v) = e^{-a_1 |v|^2 +c_1}\) and \(0<a_1<a_0\), \(c_1\),
\(\rho_0\), \(C_0\), \(C_1\) are constants. Then for any 
\(0<a<a_1\), for any \(t\in [0,T]\)
%\begeq\label{upp-point}
\[ 
f(x,v,t) \le M(v), \quad \text{for a. a. }
(x,v)\,\in\,\T^\di\times\R^\di, 
\] 
where \( M(v) = e^{-a |v|^2 +c}\), and the constant \(c\) depends on
\(a\), \(a_0\), \(c_0\), 
\(a_1\), \(c_1\), \(\rho_0\), \(C_0\) and \(C_1\) only.    
%\endeq 
\end{theorem}

\begin{remark} The regularity assumptions in 
Theorem~\ref{conditional} are not particularly 
restrictive. The precise conditions in the spatially 
inhomogeneous case are that \(f\) is a mild (renormalized)
solution of~\eqref{eq:boltz} that is dissipative in 
the sense of P.-L. Lions 
(see Definition~\ref{def:dis} in Section~\ref{sec:cp}). 
A sufficient condition that is naturally satisfied 
in the spatially-homogeneous case is 
that~\eqref{eq:ql1} holds in addition to~\eqref{eq:mild}. 
\end{remark}

The plan of the paper is as follows. In 
Section~\ref{sec:mom} we extend property (iv) 
from Theorem~\ref{alreadyknown} to the class
of Boltzmann kernels 
satisfying~\eqref{eq:asmp_b}--\eqref{eq:asmp_h}.
This involves a rather technical analysis, 
specific to the spatially-homogeneous problem,
which is based mainly on a development of the ideas 
from~\cite{Bo0,Bo,BoGaPa}. The result of 
Section~\ref{sec:mom}, however, illustrates an 
important point that the type of behavior 
described by Theorem~\ref{thm:1} is not a 
particular feature of the hard-sphere model, 
but rather a generic phenomenon that holds for a
wide class of collision kernels of ``hard'' type.   
The key step occurs 
in Section~\ref{sec:cp}: there we introduce the 
technique based on a comparison principle 
which plays a crucial role in the derivation of 
pointwise estimates. In Section~\ref{sec:wt}
we prove a weighted bound for the collision
term, based on the Carleman representation 
of the gain operator, which is used in the 
comparison argument. Finally, some classical 
results used throughout the text are recalled in 
three Appendices.
\sm

\noindent
{\bf Convention:} Throughout the text, the function \(\sign\, z\)
is defined as \(1\) for \(z>0\), \(-1\) for \(z<0\) and 
an arbitrary fixed value in \([-1,1]\) for \(z=0\). 

%The remaining of this paper is devoted to the proof of 
%Theorem~\ref{conditional}.%
%
%The proof involves two main tools. First, a novel comparison 
%principle, based on a certain monotonicity property 
%of the linear Boltzmann equation. 
%Some of the relevant estimates were 
%earlier suggested by one of the 
%authors~\cite[Chapter~2, Section~6.2]{Vi};
%Second, we need to make use of weighted estimates based on 
%the Carleman representation
%of the Boltzmann collision operator. These tools are investigated in
%Sections~\ref{secmp} and~\ref{seccarl} respectively; 
%then in Section~\ref{secccl} we briefly conclude the proof.
%
%\newpage

\section{Weighted \(L^1\) estimates of  solutions}
\label{sec:mom}

The aim of this section is the following 
result, originally due to Bobylev in the case of 
the ``hard spheres'' and Maxwell 
molecules~\cite{Bo0,Bo}.   
%
%Our assumptions on the kernel \(B\) in the Boltzmann 
%collision operator~\eqref{eq:co} are that 
%\begin{equation}
%\label{eq:asmp_b}
%B(v-v_*,\sigma) = |v-v_*|^\beta\, h(\cos\th), 
%\quad
%\cos\th = \tfrac{(v-v_*)\cdot\sigma}{|v-v_*|},  
%\end{equation}
%where \(0<\beta\le 1\) and \(h\) is a nonnegative 
%function that satisfies
%\begin{equation}
%\label{eq:asmp_h}
%h(\cos\th)\,\sin^{\alpha}\th \le C,\quad \th\in (0,\pi),
%\end{equation}
%where \(\alpha<\di-1\) and \(C\) is a constant. 
%
\begin{theorem}
\label{thm:int} 
Let \(f(v,t)\), \(v\in\R^\di\), \(t\ge0\) (\(n\ge2\)) 
be a solution of the spatially homogeneous Boltzmann 
equation~\eqref{eq:bosh} with the collision kernel 
\(B\) satisfying~\eqref{eq:asmp_b}--\eqref{eq:asmp_h}
and with the initial datum \(f_0\ge 0\) such that
\begin{equation}
\label{eq:mini}
\frac{f_0}{M_0} \in L^1(\R^\di)
\end{equation}
for a certain Maxwellian \(M_0(v)=e^{-a_0 |v|^2}\), 
where \(a_0\) is a positive constant. Then there 
exist constants \(C\), \(a>0\), such that 
\begin{equation}
\label{eq:mi_est}
\ir \frac{f(v,t)}{M(v)}\, dv \le C, \quad t\ge 0,
\end{equation}
where \(M(v)=e^{-a|v|^2}\). 
\end{theorem} 

Our approach to the problem is based on the analysis 
of the sequence of moments, 
\begin{equation}
\label{eq:mom}
m_k(t) = \ir f(v,t)\,|v|^{2k}\, dv, \quad k=0,1\dots,
\end{equation}
and particularly, of the growth of \(m_k(t)\) as \(k\to\infty\). 
The relation between the moments~\eqref{eq:mom} and 
the weighted averages~\eqref{eq:mi_est} is given by 
the formal expansion  
\begin{equation}
\label{eq:av_exp}
\ir  \frac{f(v,t)}{M(v)}\,\,dv 
\,=\, \sum_{k=0}^\infty \,\frac{m_k(t)}{k!}\,\, a^k.
\end{equation}
In view of~\eqref{eq:av_exp}, to prove Theorem~\ref{thm:int} 
it suffices to show
\begin{equation}
\label{eq:limm}
\limm\limits_{\rule{0pt}{7pt}k\to\infty} 
\,\sup\limits_{t\ge 0}\,
\frac{m_k(t)}{k!\, A^k}\, =\,0, 
\quad 
\text{for some \(A\) large enough}. 
\end{equation}

Our proof of~\eqref{eq:limm} is to a large extent a 
refinement of the 
original approach in~\cite{Bo}. One particular technical 
aspect which allows us to simplify some of the arguments
is the systematic use of the interpolation inequalities
\begin{equation}
\label{eq:mint}
\Bigl(\tfrac{m_{k_1}(t)}{m_0}\Bigr)^{\frac{1}{k_1}} 
\le 
\Bigl(\tfrac{m_{k}(t)}{m_0}\Bigr)^{\frac{1}{k}} 
\le
\Bigl(\tfrac{m_{k_2}(t)}{m_0}\Bigr)^{\frac{1}{k_2}}, 
\quad
k_1\le k \le k_2,
\end{equation} 
which follow directly from~\eqref{eq:mom} by application 
of either H\"older or Jensen's inequalities. 

It is well-known that if the kernel 
\(B(|v-v_*|,\cos\th)\) in~\eqref{eq:co} is constant 
in the first argument (the case of the Maxwell, 
or pseudo-Maxwell, particles) then the equations 
for the moments \(m_k(t)\) with integer \(k\) form 
a closed infinite system of ODE. This 
property no longer holds if the kernel \(B\)
depends on \(|v-v_*|\), and one has to work with inequalities 
instead of equations. If the kernel \(B\) has the 
homogeneity 
\(|v-v_*|^\beta\), one also generally has to consider 
the moments 
\begin{equation}
\label{eq:mos}
m_k(t)\quad\text{with} \quad 
k=j+\tfrac{\beta}{2}\, l, \quad j,\, l = 0, 1\dots
\end{equation}
Since the total mass is conserved, \(m_0(t)=m_0={\rm const}\);
we shall enumerate the rest of the moments~\eqref{eq:mos} by 
a single index \(k_n\), \(n=1,2\dots\), in the increasing
order, and introduce the notation
\[
J=\{k_n:\,n=1,2\dots\}
\] 
for the index set.
Also, introduce the normalized moments 
\begin{equation}
\label{eq:zk}
z_{k}(t) = \frac{m_{k}(t)}{\Gamma(k+b)},\quad k\ge 0,
\end{equation}
where the constant \(b>0\) will be chosen below depending 
on \(\alpha\) in \eqref{eq:asmp_h}. For \(b=1\) and 
\(k\) nonnegative integer we have \(z_{k}(t)=m_{k}(t)/k!\) 
which is the normalization appearing in~\eqref{eq:limm}. 
Also, as is easy to verify by Stirling's formula,
\begin{equation}
\label{eq:ga}
\Gamma(k+b)\sim k^{b-1}\, \Gamma(k+1), \quad k\to \infty,
\end{equation}
so the particular choice of \(b\) is irrelevant 
for~\eqref{eq:limm}. 

Given \(k=k_n\in J\) we set
\(\bar{z}^{(k)}(t)=(z_{k_1}(t),\dots,z_{k_{n-1}}(t))\), a 
vector with \(n-1\) components. 

By the assumptions on \(m_k(0)\), we have
\begin{equation}
\label{eq:zini}
z_{k}(0) \le C_0 \,q_0^{k}, \quad k\in J,
\end{equation}
for certain constants \(C_0\), \(q_0\). We shall show that the 
geometric growth of the normalized moments is preserved 
uniformly in time, due to the structure of the system of 
differential inequalities satisfied by \(z_{k}(t)\);
this will imply~\eqref{eq:limm}. 
The key step is the following.  

\begin{lemma}
\label{lem:di}
Let the sequence of nonnegative functions 
\(z_{k}\in C^1([0,\infty))\), \(k\in J\), satisfy 
\begin{equation}
\label{eq:zkine}
z'_{k}(t) \le - \,A_k \,z_{k}^{1+\frac{\beta}{2k}}(t) 
+ B_k \,F_k (\bar{z}^{(k)}(t)),\quad k\in J,\quad k\ge k_*
\end{equation}
and 
\begin{equation}
\label{eq:zklow}
z_k(t) \le C_1\, q_1^k,\quad k\in J, \quad  k < k_*, 
\end{equation}
where \(k_*>\tfrac{\beta}{2}\), \(C_1\) and \(q_1\) are 
positive constants,  \(A_k\), \(B_k\) 
are positive sequences satisfying 
\begin{equation}
\label{eq:lbk}
\frac{A_k}{B_k} \ge C_1^{1-\frac{\beta}{2k}}, \quad k\in J,\quad k \ge k_* ,
\end{equation}
and \(F_k\)  are continuous functions of their arguments 
such that
\begin{equation}
\label{eq:fbd}
F_k (\bar{z}^{(k)}) \,\le\, C^2\, q^{k+\frac{\beta}{2}},
\quad\text{whenever}\quad
z_{k} \le C q^{k}, \quad k\in J, \quad k\ge k_*.
\end{equation}
Assume that the sequence \(z_{k}(0)\) satisfies~\eqref{eq:zini}. 
Then  \(z_{k}(t) \le C q^{k}\), \(k\in J\), \(t\ge0\), 
where \(C=\max\{C_0,C_1\}\) and \(q=\max\{q_0,q_1\}\).   
\end{lemma}

\begin{proof}
Without loss of generality we can assume that 
\(C_1=C_0\) and \(q_1=q_0\). The proof will be achieved 
by induction on \(k\in J\), \(k\ge k_*\). For \(k= k_*\)
conditions~\eqref{eq:zklow} and~\eqref{eq:fbd} imply
\[
z'_k(t) \le - \,A_k \,z_{k}^{1+\frac{\beta}{2k}}(t) 
+ B_k \,C_1^2\, q_1^{k+\frac{\beta}{2}}.
\]
By a comparison argument for Bernoulli-type ordinary differential
equations (cf.~\cite{Bo}),
\begin{equation}
\label{eq:zkb}
z_k(t) \le \max\{z_k(0),z^*_k\},
\end{equation}
where \(z^*_k\) is determined from the equation
\[
\,A_k \,(z^*_{k})^{1+\frac{\beta}{2k}} 
= B_k \,C_1^2\, q_1^{k+\frac{\beta}{2}}
\]
Using condition~\eqref{eq:lbk} it is easy to verify that 
\(
z^*_{k} \le C_1\, q_1^k,
\)
which in view of~\eqref{eq:zkb} and~\eqref{eq:zini} implies
\(
z_k(t) \le C_1 q_1^k
\), \(k=k_*\). This provides the basis for the induction.
The induction step follows by repeating the same reasoning 
for any \(k>k_*\). The proof is complete. 
\end{proof}

Next, we shall verify the conditions of Lemma~\ref{lem:di} for the 
sequence of the moments corresponding to a solution of the 
Boltzmann equation. The proof of the time-regularity 
of the moments is standard; we refer the reader to
Appendix B for the details. We can also use the known 
property that the moments of every order are uniformly 
bounded in time (part (iii) of Theorem~\ref{alreadyknown}) 
to deduce~\eqref{eq:zklow}. The main difficulty is 
then to obtain the system~\eqref{eq:zkine} and to make sure
that the necessary estimates hold for the constants.

%As a consequence of Lemma~\ref{lem:mdiff} we find that 
%%for a solution of the Boltzmann equation, 
%if all 
%moments \(m_k(t)\) are finite initially then they 
%are continuously differentiable as functions of time. 
%Our next aim is to obtain a system of inequalities 
%for the normalized moments \(z_k(t)\) in the form 
%given in Lemma~\ref{lem:di}. 
%
%We follow the method of proof in~\cite{Bo}, writing
%\[
%\int_{\R^n}  \frac{f(v,t)}{M(v)}\,dv 
%= \sum_{k=0}^\infty \frac{m_k(t)}{k!}\, a^k
%\]
%where 
%\[
%m_k(t) = \ir f(v,t)\,|v|^{2k}\, dv, \quad k=0,1\dots
%\]
%For fixed \(t\) the series converges for \(a<A\) if 
%\[
%\frac{m_k(t)}{k!} \le A^k, \quad \text{for}\quad k>K_0. 
%\]
%Thus, to establish the uniform bound in 
%Theorem~\ref{thm:int} it suffices to prove
%\[
%\limm\limits_{k\to\infty} \sup\limits_{t\ge 0}\,
%\frac{m_k(t)}{k!\, A^k}\, =\,0, 
%\]
%for some \(A\) large enough. 

Let us first make some general comments about the 
time-evolution of the quantities 
\(
\ir f(v,t)\, \Psi(|v|^2)\, dv
\),
where \(\Psi:\R_+\to \R\) is a convex function. 
Multiplying equation~\eqref{eq:bosh} by \(\Psi(|v|^2)\)
and integrating with respect to \(v\) we obtain, after 
standard changes of variables,
\begin{equation}
\label{eq:evo_Psi}
\begin{split}
\frac{d}{dt} \ir f(v,t) \, \Psi (|v|^2)\,dv
 = \ir \ir f(v,t) \,f(v_*,t)\, R_\Psi (v,v_*) \, dv\, dv_*,
\end{split}
\end{equation}
where 
\[
R_\Psi(v,v_*) = 
|v-v_*|^\beta\, \big(G_\Psi(v,v_*) - L_\Psi(v,v_*)\big),
\]
\[
\rule{0pt}{22pt}
G_\Psi(v,v_*) = \frac12 \,\is \big(\Psi (|v'_*|^2) + \Psi (|v'|^2)\big)\,
h \big(\tfrac{(v-v_*)\cdot\sigma}{|v-v_*|}\big) \,d\sigma,
\]
\(v'_*\), \(v'\) are defined by~\eqref{eq:cp}, and
\begin{equation*}
L_\Psi(v,v_*) = \frac12\, \big(\Psi (|v|^2)
+ \Psi (|v_*|^2)\big).
\end{equation*}
Since the expression for $G_\Psi(v,v_*)$ is clearly the 
most complicated part of~\eqref{eq:evo_Psi} we 
look for a simpler upper bound.
This is achieved by means of the following 
estimate. 

\begin{lemma} Let \(\Psi:\R_+\to \R\) be 
convex and assume that the function 
$\bar{h}(z)=\tfrac{1}{2}(h(z)+h(-z))$ 
is nondecreasing on $(0,1)$. Then
\label{lem2}
\begin{equation*}
G_\Psi(v,v_*) \le \omega_{\di-2}\int_{-1}^1 \Psi \Big(\big( |v|^2 
+ |v_*|^2\big)\,\frac{1+z} 2\,\Big) \,
\bar{h} (z)\,(1-z^2)^{\frac{\di-3}{2}}\,dz.
\end{equation*}
\end{lemma}
\begin{proof} See~\cite[Lemma 1]{BoGaPa} for the case \(\di=3\);
the extension to general \(n\) is straightforward. 
\end{proof}

Recall that the mass \(m_0\) and the energy \(m_1\) are 
constant for the solution \(f(v,t)\). We will also use 
a lower bound for the moments of order \(\alpha\le 1\). 
\begin{lemma}[Cf. \cite{Bo} for the case \(\alpha=1\)]
\label{lem3}
The solution of~\eqref{eq:bosh} satisfies
\[
\ir f(v_*,t)\,|v-v_*|^\alpha\, dv_* \ge
c_\alpha\ir f_0(v_*)\,|v-v_*|^\alpha\, dv_*, \quad v\in \R^\di,
\]
for any \(\alpha\in(0,1]\). 
\end{lemma}
\begin{proof} By translating the solution
\(f(v_*,t)\) in the velocity space, we can reduce the proof
to the case \(v=0\). We will establish the estimates
\begin{equation}
\label{eq:ma_low}
m_\alpha (t) \ge c_\alpha\, m_\alpha(0),
\end{equation}
for \(0<\alpha\le1\). Notice that  $\Psi (z) = -z^\alpha$ 
is a convex function.
Then, by the previous computation, and using Lemma~\ref{lem2},
\[
m'_\alpha (t) \ge \ir\ir f(v,t)\, f(v_*,t)\,\,|v-v_*|^\beta\,
\Bigl( \frac{a_\alpha}2 \,\bigl(|v|^2 + |v_*|^2\bigr)^\alpha -\frac12 \,
\bigl(|v|^{2\alpha} + |v_*|^{2\alpha}\bigr)\!\Bigr)\, dv\,dv_*
\]
where $a_\alpha = 2\int_{-1}^1 (\frac{1+z}2)^\alpha 
\,\bar{b}(z)\,(1-z^2)^{\frac{\di-3}{2}}\,dz >1$.
We shall estimate the integrand above in order to obtain an 
expression involving \(m_\alpha (t)\) and similar quantities. 
For this we notice that since 
\((x+y)^\beta \le x^\beta + y^\beta\),
for \(\beta \in [0,1]\), then 
\[
|v-v_*|^\beta \le (|v| + |v_*|)^\beta \le |v|^\beta + |v_*|^\beta. 
\]
Also,
\[
|v-v_*|^\beta \ge \big|\,|v|^\beta-|v_*|^\beta\big|
\quad
\text{and}
\quad
(|v|^2+|v_*|^2)^\alpha \ge |\,|v|^{2\alpha}-|v_*|^{2\alpha}\big|. 
\]
Therefore 
\begin{equation*}
\begin{split}
&|v-v_*|^\beta\, \Bigl( \frac{a_\alpha}2\,  
\bigl(|v|^2 + |v_*|^2\bigr)^\alpha 
- \frac12\, \bigl(|v|^{2\alpha} + |v_*|^{2\alpha}\bigr)\Bigr)\\
\noalign{\vskip6pt}
&\qquad \ge \frac{a_\alpha} 2\, \bigl(|v|^\beta - |v_*|^\beta\bigr) 
\bigl(|v|^{2\alpha} - |v_*|^{2\alpha}\bigr) 
- \frac12\,  \bigl(|v|^\beta + |v_*|^\beta) 
\bigl(|v|^{2\alpha} + |v_*|^{2\alpha}\bigr)\cr
\noalign{\vskip6pt}
&\qquad = \frac{a_\alpha -1}2\,  
\bigl(|v|^{\beta +2\alpha} + |v_*|^{\beta +2\alpha}\bigr) 
- \frac{a_\alpha +1}2 \, 
\bigl(|v|^\beta |v_*|^{2\alpha} + |v|^{2\alpha} |v_*|^\beta\bigr)
\end{split}
\end{equation*}
and we obtain 
$$m'_\alpha(t) \,\ge\, (a_\alpha -1)\, m_0 \,m_{\alpha +\frac{\beta}2} (t) 
- (a_\alpha+1)\, m_{\frac{\beta}2} 
(t) \,m_\alpha (t)\  .$$
In the particular case $\beta =1$ we have
$$m'_{\frac12} (t) \ge (a_{\frac12}-1)\, m_0 m_1 
- (a_{\frac12} +1)\, m_{\frac12}^2(t),
$$
(\(m_0\) and \(m_1\) are constants, by the conservation of mass and
energy).
Therefore,
\begin{equation*}
\begin{split}
m_{{\frac12}} (t) 
\ge \min \biggl\{ m_{{\frac12}} (0), 
\biggl( \frac{a_{{\frac12}}-1}{a_{{\frac12}} +1} \,m_0 \,m_1\biggr)^{{\frac12}}
\biggr\}  
 \ge \min \biggl\{ 1, 
\biggl( \frac{a_{{\frac12}} -1}{a_{{\frac12}}+1}\biggr)^{\frac12}\biggr\} \,
m_{{\frac12}} (0)\ ,
\end{split}
\end{equation*}
since 
\(
m_0 m_1 \ge m_{{\frac12}} (0)^2
\). 
(This is the argument of Bobylev.) To achieve the proof 
for \(\beta<1\) we iterate this argument, applying it with 
\(\alpha=\frac{j\beta}2\), \(j=1\dots\), until
\(\frac{(j+1)\beta}{2}\ge1\). Consider first the case of 
the terminal \(j\), when
$$\alpha_0 = \frac{j\beta}2  < 1\le \frac{(j+1)\beta}2\ .$$
In that case
\begin{equation*}
\begin{split}
m'_{\alpha_0} (t) 
& \ge (a_{\alpha_0} -1)\, m_0\, m_{{\alpha_0} +\frac{\beta}2} (t)
- (a_{\alpha_0} +1)\, m_{\beta/2} (t)\, m_{\alpha_0} (t)\\
\noalign{\vskip6pt}
& \ge (a_{\alpha_0} -1)\, m_0^{2-({\alpha_0}+\frac{\beta}2)}\, 
m_1^{{\alpha_0}+\frac{\beta}2}
- (a_{\alpha_0} +1)\, m_0^{1-\frac{\beta}{2{\alpha_0}}}\, 
m_{\alpha_0}^{1+\frac{\beta}{2{\alpha_0}}} (t)
\end{split}
\end{equation*}
Therefore,
\begin{equation*}
 \begin{split}
m_{\alpha_0} (t) 
&\ge \min \biggl\{ m_{\alpha_0} (0), 
\biggl( \frac{a_{\alpha_0}-1}{a_{\alpha_0}+1} 
m_0^{(\frac1{{\alpha_0}}-1)({\alpha_0}+\frac{\beta}2)} 
m_1^{{\alpha_0}+\frac{\beta}2}\biggr)^{\frac1{1+\frac{\beta}{2{\alpha_0}}}}
\biggr\}
\\
\noalign{\vskip6pt}
&\ge \min \biggl\{ 1,
\biggl(\frac{a_{\alpha_0}-1}{a_{\alpha_0} 
+1}\biggr)^{\frac1{1+\frac{\beta}{2{\alpha_0}}}}
\biggr\} m_{\alpha_0} (0)
 = \biggl( \frac{a_{\alpha_0} -1}{a_{\alpha_0} 
+1}\biggr)^{\frac1{1+\frac{\beta}{2{\alpha_0}}}}
m_{\alpha_0} (0)\ .
\end{split}
\end{equation*}
Further, take $\alpha_1 = \alpha_0 - \frac{\beta}2 >0$.
Then
$$m'_{\alpha_1} (t) \ge (a_{\alpha_1} -1)\, m_0\, m_{\alpha_0} (t) 
- (a_{\alpha_1} +1)\, m_0^{1-\frac{\beta}{2\alpha_1}} \,
m_{\alpha_1}^{1+\frac{\beta}{2\alpha_1}} (t),$$
so 
\begin{equation*}
\begin{split}
m_{\alpha_1}(t)
&\ge \min \biggl\{ m_{\alpha_1} (0), 
\biggl( \biggl( \frac{a_{\alpha_1}-1}{a_{\alpha_1}+1}\biggr) 
m_0^{\frac{\beta}{2\alpha_1}} m_{\alpha_0} (t)\biggr)
^{\frac1{1+\frac{\beta}{2\alpha_1}}} \biggr\}\\
\noalign{\vskip6pt}
&\ge \min \biggl\{ m_{\alpha_1} (0),  
\biggl( \biggl( \frac{a_{\alpha_1}-1}{a_{\alpha_1}+1}\biggr)
\biggl( \frac{a_{\alpha_0}-1}
{a_{\alpha_0}+1}\biggr)^{\frac{\alpha_0}{\alpha_0+\frac{\beta}2}}
m_0^{\frac{\beta}{2\alpha_1}} 
m_{\alpha_0}(0)\biggr)^{\frac1{1+\frac{\beta}{2\alpha_1}}}
\biggr\}\\
\noalign{\vskip6pt}
&\ge \biggl( \frac{a_{\alpha_1}-1}
{a_{\alpha_1}+1}\biggr)^{\frac{\alpha_1}{\alpha_1+\frac{\beta}2}}
\biggl( \frac{a_{\alpha_0}-1}
{a_{\alpha_0}+1}\biggr)^{\frac{\alpha_1}{\alpha_0+\frac{\beta}2}}
m_{\alpha_1}(0)\ .
\end{split}
\end{equation*}
The rest of the proof can be obtained by induction. 
This establishes \eqref{eq:ma_low} for all \(\alpha
\in (0,1]\) and completes the proof.
\end{proof}

In the particular case \(\Psi(z)=z^k\), \(k\ge 1\), we 
obtain the following inequalities
\begin{equation}
\label{eq:gmk}
m'_k(t) \le  \ir\ir
f(v,t) f(v_*,t)\, \bar{R}_k(v,v_*)\,dv\,dv_*,
\end{equation}
where 
\begin{equation}
\label{grk}
\bar{R}_k(v,v_*) = \frac12\,
 |v-v_*|^\beta
\bigl( a_k  (|v|^2 +|v_*|^2)^k 
- |v|^{2k} - |v_*|^{2k}\bigr),
\end{equation}
and the constant \(a_k\) is defined by 
\begin{equation}
\label{eq:ak0}
a_k = \omega_{\di-2}\int_{-1}^1 
\Bigl( \frac{1+z}2\Bigr)^k\, \bar{h}(z)\,
(1-z^2)^{\frac{\di-3}{2}}\,dz, 
\end{equation}
satisfies \(a_1=1\), \(a_k\le 1\) for \(k\ge 1\) and is 
strictly decreasing with increasing \(k\). We notice 
the inequalities 
$|v-v_*|^\beta \le |v|^\beta + |v_*|^\beta$,   
\begin{equation*}
(|v|^2 + |v_*|^2)^k - |v|^{2k} - |v_*|^{2k} 
\le \sum_{i=1}^{[\frac{k+1}2]}
\binom{k}{j} (|v|^{2j} |v_*|^{2(k-j)} + |v|^{2(k-j)}
|v_*|^j),
\end{equation*}
where \([\,\cdot\,]\) denotes the integer part 
(cf.~\cite{BoGaPa}). Also, by Lemma~\ref{lem3},
\begin{equation*}
\begin{split}
\ir f(v_*,t)\, |v-v_*|^\beta \, dv_*
 \ge c_\beta \int f_0 (v_*)\, |v-v_*|^\beta \,dv_*
 \ge \nu_0 \,(1+|v|^\beta),
\end{split}
\end{equation*}
where \(\nu_0\) is a constant depending on \(\beta\) 
and \(f_0\). Using these inequalities 
in~\eqref{eq:gmk},~\eqref{grk} we obtain
\begin{equation}
\label{eq:mks}
m'_k(t) \le - (1-a_k)\, \nu_0\, m_{k+\frac{\beta}2} (t) 
+ a_k\, S_k(t)
\end{equation}
where 
\[
S_k(t) = \sum_{j=1}^{[\frac{k+1}2]} \binom{k}{j}
\bigl(m_{j+\frac{\beta}2} (t)\,m_{k-j}(t) 
+m_{k-j+\frac{\beta}2} (t) \,m_j (t)\bigr). 
\]
The crucial estimate for the sum \(S_k(t)\) is 
provided by the following Lemma. 

\begin{lemma} 
\label{lem:zk}
For \(b>0\)  fixed  set 
$z_k(t) = {m_k(t)}/{\Gamma (k+b)}$, $k\ge 1$. 
Then 
\[
S_k(t) \le C_b\,\Gamma (k+\frac{\beta}2 +2b)\, Z_k(t), 
\quad k\ge 1,
\]
where 
\begin{equation}
\label{eq:Zk}
Z_k(t) = \max_{1\le j\le [\frac{k+1}2]} 
\{z_{j+\frac{\beta}2}(t)\, z_{k-j}(t), 
z_j(t)\, z_{k-j+\frac{\beta}2}(t)\}
\end{equation}
and \(C_b\) is a constant depending on \(b\). 
\end{lemma}
\begin{proof} See~\cite[Lemma 4]{BoGaPa}. 
\end{proof}

\begin{proof}[Proof of Theorem~\ref{thm:int}]
Using the interpolation inequality 
\(m_{k+\frac{\beta}{2}}(t) \ge m_0^{-\frac{\beta}{2k}}\,
m_k(t)^{1+\frac{\beta}{2k}}\) and Lemma~\ref{lem:zk}
we derive from~\eqref{eq:mks} the inequalities
\begin{equation}
\label{eq:zk3}
z'_k(t) \le - (1-a_k) \,\nu_0\, 
m_0^{-\frac{\beta}{2k}}\, 
\Gamma (k+b)^{\frac{\beta}{2k}}\, 
z_k^{1+\frac{\beta}2} (t) 
+ a_k \,C_b \,
\frac{\Gamma(k+\frac{\beta}2 +2b)}{\Gamma (k+b)}\, Z_k(t).
\end{equation}
Notice that for \(k\in J\), \(k>1+\frac{\beta}{2}\)
the term \(Z_k(t)\) is of the form 
\(F_k(\bar{z}^{(k)}(t))\) as in Lemma~\ref{lem:di},
since the highest order of moment 
entering~\eqref{eq:Zk} is \(k-1+\frac{\beta}{2}\).  
It is also clear the the function \(F_k\) defined 
in this way is a continuous function of its 
arguments. Thus, we can identify~\eqref{eq:zk3}
with~\eqref{eq:zkine} by setting 
\begin{equation}
\label{eq:AkBk}
A_k = (1-a_k) \,\nu_0\, 
m_0^{-\frac{\beta}{2k}}\, 
\Gamma (k+b)^{\frac{\beta}{2k}},
\qquad
B_k = a_k \,C_b \,
\frac{\Gamma(k+\frac{\beta}2 +2b)}{\Gamma (k+b)}.
\end{equation}

We would like to apply Lemma~\ref{lem:di} to the 
sequence of functions \(z_k(t)\). 
%Conditions~\eqref{eq:zini} and~\eqref{eq:zklow} 
%are satisfied by the assumption~\eqref{eq:mini} and 
%using the fact that the moments \(m_k(t)\)
%are uniformly bounded as functions of time for 
%every \(k\), cf.~\cite{El}. 
It  remains to verify 
that the sequences of constants \(A_k\) and \(B_k\)
appearing in~\eqref{eq:zk3} satisfy~\eqref{eq:lbk}.
To this end we show that
\begin{equation}
\label{eq:lbq}
\frac{A_k}{B_k} \ge c_0, \quad\; k>k_* ,
\end{equation}
for any \(k_*>1+\frac{\beta}{2}\) and a sufficiently 
small \(c_0\); then~\eqref{eq:lbk}
would follow by choosing
\(C_0=C_1\) small enough
and \(q_0=q_1\) large enough in~\eqref{eq:zini},~\eqref{eq:zklow}. 
Indeed, using~\eqref{eq:ga},
\begin{equation}
\label{eq:Gammac}
\Gamma (k+b)^{-\frac{\beta}{2k}} \sim k^{\beta/2}
\quad
\text{and} 
\quad
\frac{\Gamma (k+\frac{\beta}2 +2b)}{\Gamma(k+b)} 
\sim k^{\frac{\beta}2+b}, 
\quad 
k\to\infty. 
\end{equation}
To estimate the constant \(a_k\) in \eqref{eq:zk3} we
recall that \(\bar{b}(z)\le C\,(1-z^2)^{-\alpha}\), 
\(\alpha<\di-1\) and setting in \eqref{eq:ak0}, 
\(s=\frac{z+1}{2}\), \(\eps={n-1-\alpha}>0\),  
we have 
\begin{equation}
\label{eq:aksymp}
\begin{split}
a_k =  C \,2^{-1+\eps} 
\int_{0}^{1} s^{k-1+\frac{\eps}{2}} (1-s)^{-1+\frac{\eps}{2}}\, ds
=  C\, 2^{-1+\eps} B(k+\tfrac{\eps}{2},\tfrac{\eps}{2}) 
\\
= 
 C\, 2^{-1+\eps}\, 
\frac{\Gamma(k+\frac{\eps}{2})\,\Gamma(\frac{\eps}{2})}
{\Gamma(k+\eps)} \asymp k^{-\frac{\eps}{2}}, \quad k\to\infty.
\end{split}
\end{equation}
We fix \(0<b< \eps/2\); the corresponding constants 
\(A_k\), \(B_k\) satisfy the inequalities  
\begin{equation}
\label{eq:ci}
A_k \ge \bar{A} \,k^{\frac{\beta}{2}}, 
\quad
B_k \le \bar{B} \,k^{\frac{\beta}{2}+b-\frac{\eps}{2}},  
 \quad k\ge  k_*,
\end{equation}
where \(k_*>1+\frac{\beta}{2}\), and \(\bar{A}\) 
and \(\bar{B}\) are absolute constants which can 
be estimated based on~\eqref{eq:AkBk} and the 
asymptotic relations~\eqref{eq:Gammac} 
and~\eqref{eq:aksymp}. From~\eqref{eq:ci} we 
obtain~\eqref{eq:lbq} by choosing
\(c_0 = \bar{A}{\bar{B}}^{-1} k_*^{\frac{\eps}{2}-b}\). 

We conclude the proof of Theorem~\ref{thm:int} by applying 
Lemma~\ref{lem:di}. 
\end{proof}

\section{Comparison principle for the Boltzmann equation}
\label{sec:cp}

In this section we discuss the important technique of
comparison that will allow us to obtain pointwise 
estimates of the solutions. The crucial property 
of the Boltzmann equation used here is a certain 
monotonicity of a {\em linear} Boltzmann semigroup. 
The argument is roughly as follows: if \(f\) is a 
solution of~\eqref{eq:boltz}, \(f|_{t=0} = f_0\), 
and \(g\) is sufficiently regular and satisfies
\begin{equation}
\label{eq:bl}
(\partial_t + v \cdot \nabla_x)\, g \ge  Q(f,g), 
\quad g|_{t=0} = g_0,
\end{equation}
then \(u=f-g\) is a solution of 
\begin{equation}
\label{eq:bd}
(\partial_t + v \cdot \nabla_x)\, u \le Q(f,u), 
\quad u|_{t=0} = u_0,
\end{equation}
where \(u_0=f_0-g_0\). We will show that if 
\(f\) is nonnegative (and satisfies certain 
minimal regularity conditions), then solutions 
of \eqref{eq:bd} satisfy the order-preserving property, 
\begin{equation}
\label{eq:mon}
\text{if} \quad u_0\le 0 \quad \text{then}\quad u\le 0
\end{equation}
(zero on the right-hand side can be replaced by any other 
solution \(\tilde{u}\) of~\eqref{eq:bd}).
This translates into the following estimate (comparison 
principle):
\begin{equation}
\label{eq:fes}
\text{if} \quad f_0\le g_0 \quad \text{and\;\; \(g\) 
satisfies~\eqref{eq:bl},\;\; then}\quad f\le g.
\end{equation}
By reversing all inequalities we obtain a similar comparison 
principle that yields lower bounds of solutions. 
 
Of course, the above scheme has to be implemented with 
suitable modifications. For instance, since {\em apriori} 
only limited information about \(f\) is available we will 
require that \(g\) satisfies~\eqref{eq:bl}
for {\em a class} of functions \(f\) (defined by the 
available { apriori} estimates). 
Another important refinement is to apply the 
estimate~\eqref{eq:fes} {\em locally} (in the case of
Theorem~\ref{conditional}, to a ``high-velocity tail'' 
\(\{|v|\ge R\}\)) since global bounds in all of the 
\((v,t)\)-space cannot be generally obtained by this 
technique. 
%Moreover, in the 
%spatially-homogeneous problem we are not able to 
%obtain inequality~\eqref{eq:bl} in all of the 
%\((v,t)\)-space, and 
%therefore we introduce a ``localized'' version 
%of the comparison principle, in which~\eqref{eq:fes} 
%is used for a ``high-velocity tail'' \(\{|v|\ge R\}\), 
%and a different argument is used to obtain the bounds 
%in the rest of the space. 
We refer to 
Proposition~\ref{prop:opp_loc} and the proof of 
Theorem~\ref{conditional} given below for the 
necessary details. In Theorem~\ref{thm:comp} we 
will give a rigorous statement of~\eqref{eq:fes}
in application to a  general class of weak 
solutions of~\eqref{eq:boltz} in the sense of DiPerna
and Lions~\cite{DPLi,Li}. 

The basic approach leading to applications of~\eqref{eq:fes} 
originated in the work by one of the 
authors~\cite[Sec.\,6.2]{Vi} in the context of lower 
bounds for the spatially-homogeneous equation without 
angular cutoff. It was also used to obtain lower 
bounds for solutions in a model describing inelastic 
collisions~\cite{GaPaVi}. 
Compared to these earlier versions we do not require 
in~\eqref{eq:fes} any differentiability in the 
\(v\)-variable, and we make more precise the 
minimal regularity conditions on \(f\).  
It is interesting to compare the present technique
with other methods based on monotonicity applied 
to the Boltzmann equation, in particular the one by 
Kaniel and Shinbrot~\cite{KaSh} (see also~\cite{IlSh,Go})
and the pointwise estimates by Vedenjapin~\cite{Ve}
(the result in the latter paper follows from our 
approach using \(g=e^{C(1+t)}\)). The monotonicity
property expressed by~\eqref{eq:mon} has also an 
important relation to the concept of  dissipative 
solutions introduced by P.-L. Lions~\cite{Li}.  
%
%It is based on a comparison principle which in turn uses
%a certain order-preserving property of the flow of the 
%{\em linear} Boltzmann equation. The underlying principle 
%is quite simple and general, and uses the dissipative 
%property of the collision term which has other important
%implications. 

We first explain the way to obtain~\eqref{eq:mon}. The 
bilinear form in~\eqref{eq:bl},~\eqref{eq:bd} is defined by 
\begin{equation}
\label{eq:bq}
Q(f,u)(x,v,t) 
= \ir \is (f'_* u' - f_* u)\, B(v-v_*,\sigma)
\, d\sigma\, dv_*, 
\end{equation}
where as usual, \(f'_* = f(x,v'_*,t)\), \(u' = u(x,v',t)\), 
\(f_* =f(x,v_*,t)\), \(u=u(x,v,t)\). At this point we 
do not need to assume the kernel 
\(B\) to satisfy~\eqref{eq:asmp_b}--\eqref{eq:asmp_h}; 
the argument goes through for a more general class of 
kernels with the usual symmetries, as described in~\cite{DPLi}, 
for instance.
%we can take a more general class of kernels with the 
%usual symmetries, as described in~\cite{DPLi}, for 
%instance. 

To illustrate the general principle, consider first 
the case of equality
in~\eqref{eq:bd}. Given \(T>0\) we fix the function 
\(f: \XD\times\R^\di\times\TD\to\R_+\), which we 
assume to be smooth in \((x,t)\), 
bounded and rapidly decaying for \(|v|\) large. 
We also assume that for every \(u_0\in D\subseteq 
L^1(\XD\times\R^\di)\) the initial-value problem
\begin{equation}
\label{eq:bleq}
(\partial_t + v \cdot \nabla_x)\, u = Q(f,u), 
\quad u|_{t=0} = u_0,
\end{equation}
has a unique solution \(u\in C([0,T];L^1(\XD\times
\R^\di))\), with enough regularity so that 
\begin{equation}
\label{eq:q_int}
Q^+(f,|u|), \, Q^-(f,|u|) \in L^1(\XD\times
\R^\di\times\TD).  
\end{equation} 
Thus, we have a well-defined  
flow map (or a semigroup)
\[
\Phi_t: \, D\ni u_0 \mapsto u(t,\cdot,\cdot)\in 
L^1(\XD\times\R^\di), \quad t\in \TD. 
\]
%Conditions on the regularity of \(f\) and \(u\)
%will be relaxed significantly later on; this will 
%be particularly 
%relevant from the point of view of applications to 
%the general weak solutions of the nonlinear spatially
%inhomogeneous problem. 

The map \(\Phi_t\) can be seen to satisfy the 
following nonexpansive property: for any 
\(u_0,\;{\tilde{u}_0}\in D\),
\begin{equation}
\label{eq:nonexp}
\ir\ir\big| \Phi_t(u_0) 
-  \Phi_t({\tilde{u}_0})\big| \,dv\,dx
\,\le\, \ir\ir| u_0 -  {\tilde{u}_0}| 
\,dv\,dx, \quad t\in \TD. 
\end{equation}
Indeed, set \(w=\Phi_t(u_0) -
\Phi_t({\tilde{u}_0})\); then 
\[
(\partial_t + v\cdot \nabla_x) w = Q(f,w) 
\quad\text{on}\quad \XD\times\R^\di\times (0,T)
\]
in the sense of distributions, and 
\(Q(f,w)\in L^1\) by our assumptions. By a 
standard argument, \(\forall\, t\in\TD\), 
for \mbox{a. a.} \((x,v)\) the function \(w^\sharp: 
s\mapsto w(x-(t-s)v,v,s)\), \(s\in \TD\), is absolutely 
continuous, and we can apply the chain rule (see Appendix A) 
to obtain  
\begin{equation}
\label{eq:gmild}
\frac{d}{ds}\, |w^\sharp| \,=\,
Q(f,w)^\sharp\,\,\sign \,w^\sharp, \quad s\in (0,T),
\end{equation}
where 
\(Q(f,w)^\sharp\) is defined similarly to \(w^\sharp\). 
Integrating with respect to \(s\in(0,t)\) and 
\((x,v)\in \XD\times\R^\di\) we obtain, 
after standard changes of variables,
\[
\begin{split}
%& 
\ir\ir |w(x,v,t)|\,dv\,dx 
= \ir\ir |w_0|\,dv\,dx
%\\& 
+ \int_0^t \ir\ir Q(f,w)\,
\,\sign \,w\, dv\,dx \,ds
%\,\le\, 0,
\end{split}
\]
where \(w_0 =  u_0 -  {\tilde{u}_0}\). 
We further notice that the bilinear 
collision term~\eqref{eq:bq} satisfies
\begin{equation}
\label{eq:qdiss}
\ir Q(f,u)\,\, \sign \,u\, dv \le 0, 
\end{equation}
for every \(f\ge 0\) and every \(u\) so 
that \(Q^+(f,|u|)\),  \(Q^-(f,|u|)\in L^1\). 
This follows immediately 
from the weak form
\[
\ir Q(f,u)\,\, \sign \,u\, dv = \ir \ir \is f_* u \,
(\sign \,u' - \sign \,u) \,B\, 
d\sigma\, dv_* dv
\] 
by noticing that
\(
u \,
(\sign\, u' - \sign\, u) \le 0 
\).

The same approach can be followed to obtain~\eqref{eq:mon}. 
Indeed, we have by~\eqref{eq:qdiss} and the mass conservation
\begin{equation}
\label{eq:uplus}
\ir Q(f,u)\, \tfrac{1}{2}(\sign\, u+1) \,dv\, \le 0, 
\end{equation}
where \(\tfrac{1}{2}(\sign\, u+1)\) is the \mbox{a. e.} derivative 
of the Lipschitz-continuous function \(u_+=\max\{u,0\}\). 
We then have 
\[
\frac{d}{ds}\, u_+^\sharp \,=\,
Q(f,u)^\sharp\,\,\tfrac{1}{2}(\sign\, u+1)^\sharp, \quad s\in (0,T),
\] 
and the integration yields 
\[
\ir\ir u_+(x,v,t)\,dv\,dx \le \ir\ir {u_0}_+\,dv\,dx, \quad t\in [0,T],
\]
which implies~\eqref{eq:mon} for a. a. \((x,v)\). 

\begin{remark} Relation~\eqref{eq:mon} can be
restated as the order-preserving property of \(\Phi_t\): 
\begin{equation}
\label{eq:opp}
\forall\; u_0,\, \tilde{u}_0\in D, \;\;
u_0\le \tilde{u}_0 \;\;\text{implies}\;\; \Phi_t(u_0)\le
\Phi_t(\tilde{u}_0), \quad t\in [0,T] . 
\end{equation} 
In fact, the equivalence of~\eqref{eq:opp} and~\eqref{eq:nonexp}
follows from a general principle applied to (nonlinear) 
maps that preserve integral, as described by Crandall and 
Tartar~\cite{CrTa}. Inequality~\eqref{eq:mon} (or~\eqref{eq:opp}) 
can then be seen as a consequence of the results in \cite{CrTa}, the 
preservation of the mass \(\ir\ir f\, dv\,dx\) along solutions 
of \eqref{eq:bd}, and \eqref{eq:nonexp}.
\end{remark}

The following localized version of the order-preserving 
property will be useful for the comparison argument. 

\begin{proposition}
\label{prop:opp_loc}
Let \(\;f,u\in C(\TD; L^1(\XD\times\R^\di))\;\) satisfy
\[
f\ge 0; \quad\partial_t u + v \cdot \nabla_x u, \; 
Q^+(f,u), \;Q^-(f,u) \in
L^1; \quad u|_{t=0}=u_0\le 0,
\]
and assume that for a certain (measurable) set
\(U\subseteq \XD\times\R^\di\times(0,T)\),
\[
\partial_t u + v\cdot \nabla_x u - Q(f,u) 
\,\le\, 0\quad\text{on}\quad U,
\] 
and
\[
u\le 0\quad\text{on}\quad U^c :=
\big(\XD\times\R^\di\times(0,T)\big)\setminus U.  
\]
Then \(u(t,\cdot,\cdot)\le0\;\) a. e. on \(\XD\times\R^\di\),
for every \(t\in \TD\). 
\end{proposition}
\begin{proof}
Let \(D(u) = \partial_t u + v \cdot \nabla_x u\). We obtain by arguing
as above,
\[
\begin{split}
\ir\ir u_+(x,v,t) \, dv\,dx%\,\Big|_{s=t} 
& \,-\, \ir\ir u_+(x,v,0) \, dv\,dx%\,\Big|_{s=0} 
\\
& =\, \int_0^t \ir \ir D(u)\,\tfrac{1}{2}(\sign\, u +1)\,dx\, dv\,
ds. 
\end{split}
\]
We have \(u_+|_{t=0}=0\); also \(\tfrac{1}{2}(\sign\, u +1) = 0\) 
whenever \(u< 0\) and  
\(D(u) = 0\) outside of a set of zero 
measure in  \(\{u = 0\}\). Therefore, 
setting \(U_t = \{(x,v,s)\in U: s\le t\}\) we have 
\[
\begin{split}
& \ir\ir u_+(x,v,t) \, dv\,dx = 
\iiint_{\,U_t} D(u)\,dx\, dv\,ds 
\\
& \le \,
\iiint_{\,U_t} Q(f,u)\,dx\, dv\,ds 
 \,=\, \int_0^t \ir \ir Q(f,u)\,\tfrac{1}{2}(\sign u +1)\,dx\, dv\,
ds\le 0,
\end{split}
\]
for every \(t\in[0,T]\),
where we used the dissipative property~\eqref{eq:uplus}. 
This shows that \(u(t,\cdot,\cdot)\le 0\) almost everywhere. 
\end{proof}

Proposition~\ref{prop:opp_loc} is sufficient to formulate 
the comparison principle in the generality required for 
Theorem~\ref{thm:1}. We will, however, give a  more general 
statement that applies to weak solutions in the spatially
inhomogeneous case. In the definition of weak solutions one 
has to account for the fact that the bound
\[
Q(f) \in L^1_{\rm loc} (\T^\di\times\R^\di\times(0,+\infty))
\]
is generally not available, and one has to define solutions 
in a sense that is weaker than distributional. The simplest
way to state the definition is to require that \(f\ge 0\),
\(f\in C([0,T];L^1_{xv})\), 
\(Q^\pm(f)/(1+f)\in 
L^1_{\rm loc}\) and the renormalized 
form
\[
(\partial_t + v\cdot \nabla_x) \, \log(1+f) = Q(f)/(1+f)
\]
holds in the sense of distributions, cf.~\cite{DPLi}. 
Such solutions are known as renormalized. This concept 
can be further refined as follows, cf.~\cite{Li}. 
%
%Proposition~\ref{prop:opp_loc} is sufficient to formulate 
%the comparison principle in the generality required for 
%Theorem~\ref{thm:1}. However, the regularity assumptions, 
%particularly the integrability  conditions for the collision 
%terms, can be relaxed to give a more general statement 
%that is applicable to a wide class of weak solutions 
%of~\eqref{eq:boltz}. We first give a definition of weak 
%solutions. 
%
%By analogy with the definition of mild solutions 
%of~\eqref{eq:bosh} given in the introduction we say 
%that \(f\in C(\TD;L^1(\XD\times\R^\di))\)
%is a mild solution of~\eqref{eq:boltz} if \(\forall\,t\), 
%a. a. \((x,v)\) the function 
%\(Q(f)^\sharp: s\mapsto Q(f)(x-(t-s)v,v,s)\)
%is integrable over \((0,T)\) and satisfies 
%\[
%f^\sharp(t_2) - f^\sharp(t_1) = \int_{t_1}^{}
%\]
%
%We will next use the order-preserving property to obtain
%a comparison principle for weak solutions of the nonlinear 
%Boltzmann equation~\eqref{eq:boltz}. 
%%\begin{equation}
%%\label{eq:bsih}
%%\partial_t f + v\cdot\nabla_x f = Q(f). 
%%\end{equation}
%To be rigorous, we first need to give a definition of
%the concept of solutions.
%
\begin{definition}
\label{def:dis}
We say that a renormalized solution \(f\) is dissipative
if \(f|v|^2\in L^\infty([0,T]; L^1_{xv})\) and 
for every sufficiently regular function
\(g:\XD\times\R^\di\times\TD\to\R\),
\begin{equation}
\label{eq:dsol}
\partial_t\! \ir |f-g|\, dv + {\rm div}_x\! \ir |f-g|\,v\, dv 
\le \ir \big(Q(f,g)-D(g)\big)\, \sign(f-g)\, dv,
\end{equation}
in the sense of distributions, where 
\(D(g)=(\partial_t + v\cdot\nabla_x) g\), and \(\sign(0)\)
is assigned an arbitrary value in \([-1,1]\). 
\end{definition}

\begin{remark}
In the above definition ``sufficiently regular''  precisely 
means that \(g\in C([0,T];L^1_{xv})\), \(g|v|^2\in 
L^\infty_t(L^1_{xv})\), 
\(D(g)\in L^1_{xvt}\) and that for any 
\(f\in C(\TD;L^1_{xv})\) such that \(f|v|^2\in L^\infty_t(L^1_{xv})\), 
\(Q^+(f,|g|)\),
\(Q^-(f,|g|)\in L^1_{xvt}\) (these conditions can 
be made more explicit, see~\cite{Li} for details).
\end{remark}

The formal motivation for the definition of dissipative 
solutions is clear: the right-hand side of the Boltzmann 
equation can be written as 
\[
Q(f) = Q(f,f-g) + Q(f,g), 
\]
so we have 
\[
(\partial_t+v\cdot\nabla_x)(f-g) = Q(f,f-g)+ Q(f,g)-D(g). 
\]
Multiplying the above equation by \(\sign(f-g)\) and using 
relation \eqref{eq:qdiss} (note that \(f\ge 0\)) we see that 
every sufficiently regular solution of~\eqref{eq:boltz} should 
satisfy~\eqref{eq:dsol}. 

Dissipative solutions are known to exist globally in time, 
for a quite general class of initial data.  
In fact, in~\cite{Li} Lions established a large class of
``dissipation inequalities'' similar to~\eqref{eq:dsol}
that hold for renormalized solutions
of~\eqref{eq:boltz}. Such solutions can also be 
constructed so that the local
mass conservation law, 
\begin{equation}
\label{eq:lmc}
\partial_t\! \ir f\, dv + {\rm div}_x\! \ir f\,v\, dv 
= 0, 
\end{equation}
holds in the sense of distributions. However they need not
generally satisfy the conditions 
\(Q^+(f),\; Q^-(f)\in L^1_{\rm loc}\).

Using the order-preserving property of 
Proposition~\ref{prop:opp_loc} we establish the
following comparison principle for dissipative solutions 
of the nonlinear Boltzmann equation. 

\begin{theorem} 
\label{thm:comp}
Let \(f\in C(\TD;L^1(\XD\times\R^\di))\) be a dissipative 
solution of \eqref{eq:boltz} and let \(g\) be a sufficiently 
regular function, such that 
\(
f|_{t=0} \le g|_{t=0},
\)
\[
\partial_t g + v\cdot \nabla_x g - Q(f,g) \ge 0
\;\;\text{on}\;\; U
\]
and 
\(
f \le g 
\;\;\text{on}\;\; U^c
\),
where \(U\) is a measurable subset of 
\(\XD\times\R^\di\times \TD\). Then \(f\le g\) almost 
everywhere on \(\XD\times\R^\di\), for every \(t\in \TD\). 
\end{theorem}

\begin{remark}
It is natural to call \(g\) a (localized) upper barrier. 
By reversing all inequalities in the above formulation 
one can also obtain a similar comparison principle for 
the lower barrier. 
\end{remark}

\begin{proof}[{\bf Proof.}] We use the notation
\(D(g)=\partial_t g + v \cdot\nabla_x g\), so that 
\[
\partial_t\! \ir g\, dv + {\rm div}_x\! \ir g\,v\, dv 
= \ir D(g)\, dv,
\]
in the sense of distributions. Using the mass 
conservation~\eqref{eq:lmc} and the identity
\[
(f-g)_+ = \tfrac{1}{2}\,\big(|f-g|+(f-g)\big)
\]
we obtain, by combining the above relations 
with~\eqref{eq:dsol}, 
\[
\begin{split}
 \partial_t\! \ir  & (f-g)_+\, dv + {\rm div}_x\! \ir (f-g)_+\,v\, dv 
\\
& \le \frac{1}{2}\ir \big(Q(f,g)-D(g)\big)\, \sign(f-g)\, dv
-\frac{1}{2}\ir D(g)\, dv.
\end{split}
\]  
Since \(Q^\pm(f,|g|)\) are integrable, we have \(\ir Q(f,g)\, dv =0\),
a. e. \((x,t)\),
and therefore,  
\begin{equation}
\label{eq:pp}
\begin{split}
 \partial_t\! \ir  & (f-g)_+\, dv + {\rm div}_x\! \ir (f-g)_+\,v\, dv 
\\
& \le \ir \big(Q(f,g)-D(g)\big)\, \tfrac{1}{2}(\sign(f-g)+1)\, dv.
\end{split}
\end{equation}
We can choose \(\sign (0)=-1\) in~\eqref{eq:pp} to avoid 
estimating the integral over the set \(\{f=g\}\). 
Since \((f-g)_+ \,v \in L^1(\XD\times\R^\di\times \TD)\) we can
integrate over \(x\) and \(t\) to obtain 
\begin{equation}
\label{eq:fgplus}
\begin{split}
\ir\ir (f-g)_+ (x,v,t)\, dv\, dx 
& \le \ir\ir (f-g)_+(x,v,0)\, dv\, dx
\\
& +\, \iiint_{\,U_t} \,\big(Q(f,g)-D(g)\big)\, dx\,dv\, ds \le 0, 
\end{split}
\end{equation}
where \(U_t = \{(x,v,s)\in U : s\le t\}\) 
and 
we used that \(\frac{1}{2}(\sign(f-g)+1)\) vanishes 
for \(f\le g\) and that \(Q(f,g)-D(g) \le 0\) 
on \(U_t\). The inequality in \eqref{eq:fgplus} implies that 
\(f\le g\), a. e. \((x,v)\in \XD\times\R^\di\), 
for every \(t\in \TD\). 
\end{proof}

Theorem~\ref{thm:comp} is a crucial ingredient in the proof of 
Theorem~\ref{conditional}, which we give below. 
\begin{proof}[Proof of Theorem~\ref{conditional}.]
To apply Theorem~\ref{thm:comp} we set 
\(U=\{(x,v,t): |v|> R\}\), where \(R\) will be chosen 
large enough, and \(g(x,v,t)=M(v)\), where 
\(M(v)=e^{-a|v|^2+c}\), \(0<a<a_1\) is fixed 
and \(c>c_0\) will be chosen sufficiently large, 
depending on \(R\). To prove that \(g\) can be used 
as a barrier for the solution on \(U\) we need to 
verify the inequality 
\begin{equation}
\label{eq:barr}
Q^+(f,g)(x,v,t) \le Q^-(f,g)(x,v,t), 
\quad (x,t)\in \T^\di\times[0,T], \quad |v| >R.  
\end{equation}
First notice that, 
by elementary inequalities, 
\[
\begin{split}
Q^-(f,g) (x,v,t)
= M(v)\,\ir f(x,v_*,t)\,|v-v_*|^\beta\, dv_*
\\
\ge M(v)\,\Big(\rho_0 |v|^\beta - 
\ir f(x,v_*,t)\,|v_*|^\beta\, dv_*\Big), 
\end{split}
\]
where \(\rho_0\) is the constant in~\eqref{eq:rho0}. 
The last term can be controlled using the estimate for the
integral of \(f/M_1\) from~\eqref{eq:fbds} as follows, 
\[
\ir f(x,v_*,t)\,|v_*|^\beta\, dv_* 
\le L \ir \frac{f(x,v_*,t)}{M_1(v_*)}\, dv_*
\le L\, C_1, 
\]
where \(L=\max\limits_{y\ge0}\, y^\beta \,e^{-a_1 y^2 + c_1}\). 
Thus, we have
\[
Q^-(f,g) (x,v,t)
\ge  M(v)\,\big(\rho_0 |v|^\beta - L\, C_1\big).  
\]
The control of the ``gain'' term is more technical; 
we establish below in Lemma~\ref{lem:w_est} the estimate
\begin{equation}
\label{eq:gainc}
Q^+(f,g) (x,v,t) \le C\,(1+|v|^{\beta-\eps})\, M(v),
\end{equation}
where \(\eps =\min\{\beta,n-1-\alpha\}>0\). This 
implies that~\eqref{eq:barr} holds if we set \(R\) 
to be the largest root of the equation
\[
C+LC_1+Cy^{\beta-\eps}-\rho_0 y^\beta =0. 
\]
Finally, we take \(c=aR^2+\log C_0\), where \(C_0\)
is the constant in~\eqref{eq:fbds}; then it is easy 
to verify that 
\begin{equation}
\label{eq:inball}
f(x,v,t) \le C_0\le  M(v), 
\quad (x,t)\in \T^\di\times[0,T], \quad |v| \le R.  
\end{equation}
The conditions \(0<a<a_1<a_0\) and  \(c\ge c_0\) guarantee 
that we have \(f(x,v,0)\le M(v)\). Together with the 
inequalities~\eqref{eq:barr} and~\eqref{eq:inball}
this allows us to use Theorem~\ref{thm:comp} 
to conclude.  
\end{proof}
%
%We further wish to apply the comparison principle 
%to the spatially homogeneous problem (i. e. \(f=f(v,t)\)
%to establish pointwise estimates for the ``tails''
%(i. e. for \(f(v,t)\), \(|v|>R\), where \(R\) is chosen 
%suitably large). 
%
%We choose the barrier function \(g=g(v)\) to be 
%time-independent 
%(for the time-uniform bounds). The inequality for 
%an upper barrier becomes simply \(Q(f,g)\le 0\), or
%\[
%Q^+(f,g)(v,t) \le Q^-(f,g)(v,t), \quad |v|\ge R, \quad t>0.  
%\]
%Clearly, an inequality in the form above cannot 
%hold for all \(v\), since \(\ir Q(f,g)\, dv =0\). 
%Remarkably, we discover that the inequality does 
%hold if we take \(R\) sufficiently large, and if \(g\)
%is a  Maxwellian function, with sufficiently 
%large density and temperature parameters. This 
%requires a certain novel weighted estimate of the 
%linear ``gain'' term which we show in the next section. 

\section{A weighted estimate for the ``gain'' operator}
\label{sec:wt}

To complete the proof of Theorem~\ref{conditional}
we prove the following weighted estimate of the 
linear ``gain'' operator. The main technique is 
based on Carleman's form of the ``gain'' term
(see Appendix C).  

\begin{lemma}
\label{lem:w_est}
Let \(B:\R^\di\times S^\dii \to \R^+\), \(n\ge2\), 
be a measurable function 
that satisfies
\[
B(u,\sigma) \,\le \,C\, (1+|u|^\beta)\,\frac{1}{|\sin \th|^\alpha}\, 
1_{\{\cos\th
  \ge 0\}}, \quad \cos\th = \tfrac{u\cdot \sigma}{|u|},
\]
where \(\beta> 0\) and \(\alpha<n-1\). Define
\[
Q^+(f,g)(v) \,=\, \ir \is f'_*\, g' \, B(v-v_*,\sigma) \, d\sigma \,
dv_*,
\]
and set \(M(v) = e^{-a |v|^2}\), \(\,a>0\);
\(\,w_\eps(v) = 1+|v|^{\beta-\eps}\), 
where \(\eps = \min\{\beta,n-1-\alpha\}> 0\). 
Then
\begin{equation}
\label{eq:we}
\Big\|\,\frac{Q^+(f,M)}{w_\eps\, M}\,\Big\|_{L^\infty(\R^\di)} 
\,\le\, C \,\Big\|\,\frac{f\, w_\eps}{M}\,\Big\|_{L^1(\R^\di)}, 
\end{equation}
where \(C\) is an explicitly computable constant depending 
on \(n\), \(\alpha\), \(\beta\) and \(a\). 
\end{lemma}
\begin{remark} For \(B\) satisfying the estimate with 
\(\alpha=0\) (for example, the kernel \(\bar{B}\) for 
hard spheres in three dimensions) we have \(\eps=\beta\) 
for all \(\beta\le \di-1\) and the weight \(w_\eps(v)\) is constant. 
The estimate of the Lemma then takes
a particularly simple form, 
\[
\Big\|\,\frac{Q^+(f,M)}{M}\,\Big\|_{L^\infty_v} 
\,\le\, C \,\Big\|\,\frac{f}{M}\,\Big\|_{L^1_v}. 
\]
For the quadratic ``gain'' term this implies the 
estimate 
\[
\Big\|\,\frac{Q^+(f)}{M}\,\Big\|_{L^\infty_v} 
\,\le\, C \,\Big\|\,\frac{f}{M}\,\Big\|_{L^\infty_v}\,
\Big\|\,\frac{f}{M}\,\Big\|_{L^1_v}. 
\]
\end{remark}
\begin{proof} By the Carleman representation formula, 
\[
Q^+(f,M)(v) = 2^{\di-1} \ir \frac{f(v'_*)}{|v-v'_*|} 
\int_{E_{vv'_*}} M(v') \,\frac{B(v-v_*,\sigma)}{|v-v_*|^{n-2}}\, 
d\pi_{v'}\,, 
\]
where \(E_{vv'_*}\) is the hyperplane 
\[
\{v'\in \R^\di \,:\, (v-v')\cdot(v-v'_*) = 0 \}, 
\]
and \(d\pi_{v'}\) denotes the usual Lebesgue measure 
on \(E_{v v'_*}\). We then have 
\begin{equation}
\label{eq:carl_int}
\frac{Q^+(f,M)(v)}{M(v)} = \ir \frac{f(v'_*)}{M(v'_*)} \, 
K(v,v'_*)\, dv'_*, 
\end{equation}
where 
\begin{equation}
\label{eq:ker_K}
K(v,v'_*) = \frac{2^{\di-1}}{|v-v'_*|} \int_{E_{vv'_*}}\!
M(v_*)\,\, \frac{B(v-v_*,\sigma)}{|v-v_*|^{n-2}} \, d\pi_{v'}, 
\end{equation}
and we used that, by the energy conservation, 
\[
\frac{M(v')\, M(v'_*)}{M(v)} = M(v_*). 
\]
Note that in \eqref{eq:ker_K} the variables 
\(v_*\) and \(\sigma\) are expressed through 
\(v\), \(v'_*\) and \(v'\) as follows, 
\[
v_* \,=\, v'_*+v'-v,\qquad
\sigma = \frac{v'-v'_*}{|v'-v'_*|}. 
\]
Now to establish the Lemma it suffices to verify the inequality
\begin{equation}
\label{eq:ker_est}
K(v,v'_*) \,\le \, C\, (1+|v-v'_*|^{\beta-\eps}). 
\end{equation}
Indeed, since 
\[
1+ |v-v'_*|^{\beta-\eps} \le
(1+|v|^{\beta-\eps})\,(1+|v'_*|^{\beta-\eps}), 
\]
then \eqref{eq:carl_int} and \eqref{eq:ker_est} imply 
\[
Q^+(f,M)(v) \,\le\, C\,(1+|v|^{\beta-\eps})\, M(v)\, \ir
\frac{f(v'_*)}{M(v'_*)}\,(1+|v'_*|^{\beta-\eps})\, dv'_*
\]
which is equivalent to~\eqref{eq:we}.  

In the remainder of the proof we will therefore verify  
\eqref{eq:ker_est}. 
Using the identity
\[
(v-v_*)\cdot (v'-v_*) \,=\, |v-v'_*|^2 - |v-v'|^2 
\]
for \(v'\in E_{vv'_*}\) and recalling that 
\(B(v-v_*,\sigma)\) vanishes for 
\((v-v_*)\cdot \sigma<0\) we see that 
the integration in \eqref{eq:ker_K} can be 
restricted to the disk 
\[
D_{vv'_*} = E_{vv'_*} \,\cap\,\{v'\in \R^\di\,:\,|v-v'_*| \le |v-v'|
\}. 
\]
We notice that for \(v'\in D_{vv'_*}\), 
\[
\bigl|\,\tan \frac{\th}{2}\,\bigr| \,=\, 
\frac{|v'_*-v_*|}{|v-v'_*|}, \qquad |\th|\le \frac{\pi}{2},
\]
where \(\th\) is the angle between the vectors 
\(v-v_*\) and \(\sigma\). This implies 
\[
\frac{1}{|\sin \th\,|} \,\le\, \frac{1}{2}\,\frac{|v-v'_*|}{|v'_*-v_*|}
\]
Thus, 
\(
\;K(v,v'_*) \le C \widetilde{K}(v,v'_*)\;\), where 
\[
\widetilde{K}(v,v'_*) = \frac{2^{\di-1-\alpha}}{|v-v'_*|^{1-\alpha}} 
\int_{D_{vv'_*}}\!
M(v_*)\,\, \frac{1+|v-v_*|^\beta}{|v-v_*|^{n-2}} \,
\frac{1}{|v'_*-v_*|^\alpha}\, d\pi_{v'}.  
\]
To estimate the above expression we consider two cases. 

{\bf Case a)} \(\;|v-v'_*|\le 1\). 
\; Since for \(v'\in D_{vv'_*}\)
\[
|v-v'_*| \le |v-v_*| \le \sqrt{2}\, |v-v'_*|
\]
we have \(1+|v-v_*|^\beta\le 1 +
2^{\beta/2}\) and 
\[
|v-v_*|^{2-n} \,\le\, |v-v'_*|^{2-n}. 
\]  
Therefore, 
\[
\widetilde{K}(v,v'_*) \,\le\, \frac{2^{\di-1-\alpha}(1 +
2^{\beta/2})}{|v-v'_*|^{n-1-\alpha}}\int_{D_{vv'_*}}\!
M(v_*)\,\,
\frac{1}{|v'_*-v_*|^\alpha}\, d\pi_{v'}. 
\]
Since \(M(v_*)\le 1\) the last integral is estimated above by  
\[
\int_{D_{vv'_*}}
\frac{1}{|v'_*-v_*|^\alpha}\, d\pi_{v'} 
= \int_{\{w\in \R^{\di-1}\,:\,|w|\le |v-v'_*|\}} 
\frac{1}{|w|^{\alpha}}\, dw = \tfrac{\omega_{\di-2}}{\di-1-\alpha}\,
|v-v'_*|^{\di-1-\alpha},
\]
if \(\di-1-\alpha>0\), i. e. \(\alpha<\di-1\). Here 
\(\omega_{\di-2}\) is the measure of the \((n-2)\)-dimensional 
unit sphere. This implies the estimate 
\[
\widetilde{K}(v,v'_*) \,\le\,\frac{2^{\di-1-\alpha}(1 +
2^{\beta/2})\,\omega_{\di-2}}{\di-1-\alpha}, \quad |v-v'_*|\le 1. 
\]

{\bf Case b)} \, \( |v-v'_*|> 1\). \; Then 
\[
1+|v-v_*|^\beta \le 2\, |v-v_*|^\beta \le  2^{1+\frac{\beta}{2}}\, |v-v'_*|^\beta,
\]
and we obtain, similarly to the previous case, 
\[
\widetilde{K}(v,v'_*) \,\le\,\frac{2^{\di-\alpha+\frac{\beta}{2}}}
{|v-v'_*|^{n-1-\alpha-\beta}}\,\int_{D_{vv'_*}}\!
M(v_*)\,\,
\frac{1}{|v'_*-v_*|^\alpha}\, d\pi_{v'}.
\]
Since \(M(v_*)\) is a radially decreasing function 
of \(v_*\in \R^\di\), and so 
is \(|v_*|^{-\alpha}\), 
\[
\begin{split}
& \int_{D_{vv'_*}}\!
M(v_*)\,{|v'_*-v_*|^{-\alpha}}\, d\pi_{v'} \le 
\int_{\R^{\di-1}} \bar{M}(w)\, |w|^{-\alpha}\, dw 
\\
& \le \int_{|w|\le 1} |w|^{-\alpha}\, dw +  \int_{\R^{\di-1}} \bar{M}(w)\, dw 
= \frac{\omega_{\di-2}}{\di-1-\alpha} 
+ \Big(\frac{\pi}{a}\Big)^{\frac{\di-1}{2}}, 
\end{split}
\]
where \(\bar{M}(w) = e^{-a|w|^2}\), \(w\in \R^{\di-1}\). 
Since \(|v-v'_*|^{\beta+\alpha-n+1}\le|v-v_*|^{\beta-\eps}\) 
this establishes the 
required estimate for Case b). 
\end{proof}

\appendix

\section*{Appendix A: Some properties of weakly 
differentiable functions}

Let \({ AC}[a,b]\) denote the class of absolutely 
continuous real-valued functions defined on an interval \([a,b]\). 
Given \(f\in { AC}[a,b]\)  we set \([c,d]=f([a,b])\)
and use the notation \({\rm Lip}[c,d]\) for the 
class of all Lipschitz continuous functions defined 
on \([c,d]\). Every function 
\(\beta\in {\rm Lip}[c,d]\)
is differentiable (in the classical sense) almost 
everywhere on \((c,d)\); we agree to extend this 
derivative to a function \(\beta'\) defined everywhere 
on \([c,d]\) by assigning arbitrary {\em finite} values 
at the points where \(\beta\) is not differentiable.  
%We will use the notations \({ AC}[a,b]\) for a set of 
%absolutely continuous functions defined on an interval 
%\([a,b]\) and \({\rm Lip}(\R)\) for the set of 
%Lipschitz continuous functions on the real line. 
%If \(\beta\in {\rm Lip}(\R)\) then it is differentiable
%almost everywhere; by redefining \(\beta'\) arbitrarily 
%on a set of measure zero if necessary we can assume 
%that it is defined at every point in \(\R\) and is finite 
%there. 
The function \(\beta'\) also coincides with the weak 
derivative of \(\beta\) almost everywhere on \((c,d)\)
The following chain rule was used in the arguments 
in Section~\ref{sec:cp}.  

\begin{proposition} Let \(f\in { AC}[a,b]\) and 
\(\beta\in {\rm Lip}[c,d]\).   
Then \(\beta\circ f\in { AC}[a,b]\) and 
\[
(\beta\circ f)' = (\beta'\circ f)\, f', 
\]
almost everywhere on \((a,b)\). 
\end{proposition}

\begin{remark} 1) The seeming ambiguity in the above 
formulation occuring since \(\beta'\circ f\)
can assume arbitrarily assigned values on a set 
of positive measure is resolved by observing 
that whenever this happens then \(f'\) vanishes, 
except on a set of measure zero (see the proof below). 
2) For the purposes of 
Section~\ref{sec:cp} we only need the chain rule 
for \(\beta(y)=|y|\) and \(\beta(y)=y_+\); these 
cases are covered in~\cite{EvGa}, and the proof for 
the case of piecewise-\(C^1\) functions \(\beta\)
can be found in~\cite{GiTr}. We include a short proof 
that applies to the general case to make 
the presentation in Section~\ref{sec:cp} self-contained.  
%By agreeing to consider \(\beta'\) everywhere 
%finite we avoid the formal problem related to interpreting 
%the right-hand side of the chain rule if \(\beta'\circ f\)
%takes the values \(\pm\infty\) on a set of positive measure. 
%As will be clear from the proof, in that case \(f'\) must 
%necessarily vanish except on a set of measure zero, and so 
%does \((\beta\circ f)'\). 
\end{remark}

\begin{proof}[{\bf Proof.}] By the definition of absolutely 
continuous functions, 
\[
\forall\, \eps>0 \;\; \exists\,\delta>0\; \text{ such that }
\;\forall\, n\in \N,\;\; \forall\, 
\{(x_j,y_j)\subseteq [a,b]: j=1,\dots, n\}, 
\]
a disjoint family, 
\[
\sum_{j=1}^{n} \,|y_j-x_j| < \delta \;\;\Rightarrow\;\;
\sum_{j=1}^{n}\, |f(y_j)-f(x_j)| < \eps. 
\]
Clearly then, since 
\[
|\,\beta(f(y_j))-\beta(f(x_j))\,| \le L\,|f(y_j)-f(x_j)|,
\]
where \(L\) is the Lipschitz constant of \(\beta\), the 
composition \(\beta\circ f\) is absolutely continuous 
on \([a,b]\). By Lebesgue's differentiation theorem, 
\(f\) and \(\beta\circ f\) are differentiable in the 
classical sense on a set with complement of measure zero 
in \((a,b)\). Pick \(x\in(a,b)\) from this set. 
We will 
consider two cases, depending on whether 
\(\beta\) is differentiable at \(f(x)\) or not. 
In the first case we have 
\[
\begin{split}
& (\beta\circ f)'(x) = \lim_{h\to 0} \frac{\beta(f(x+h))-\beta(f(x))}{h}
\\
& =\lim_{h\to 0} \frac{\beta(f(x+h))-\beta(f(x))}{f(x+h)-f(x)}
\,\lim_{h\to 0} \frac{f(x+h)-f(x)}{h}
=\beta'(f(x)) f'(x). 
\end{split}
\]
Let us further take \(A\) to be the set of \(y\) such that 
\(\beta\) is not differentiable at \(f(y)\). We claim that
\(f'(x)\) vanishes for \(x\in A\), except perhaps on a set of 
zero Lebesgue measure. Indeed, let \(B=\{y\in A: |f'(y)|>0\}\);
then 
\[
B=\mathop{\cup}\limits_{n=1}^\infty B_n, \quad B_n = 
\{y\in B:\, |f(z)-f(y)|\ge\tfrac{|z-y|}{n} 
\text{ for }|z-y|<\tfrac{1}{n}\,\}. 
\]  
We prove the claim by showing that every set \(B_n\) has 
zero measure. 

Fix an \(n\in \N\). Since \(\beta\) is 
Lipschitz, we know that 
\(f(A)\) is a set of measure zero. Given \(\eps>0\) 
we can then choose the 
intervals \(I_j\), \(j=1,\dots\), such that 
\[
f(A)\subseteq \mathop{\cup}\limits_{j=1}^\infty I_j 
\quad \text{and}\quad \sum\limits_{j=1}^{\infty}\, |I_j| <\eps. 
\]
Let \(J\) be an interval of length \(\frac{1}{n}\), and 
let \(D=B_n\cap J\), \(D_j = f^{-1}(I_j)\cap D\). Then, 
from the definition of \(B_n\), 
\(
|D_j| \le n |I_j|; 
\) therefore, \(|D|\le n\eps\) and 
\(|B_n|\le n^2|b-a|\eps\). Since \(\eps\) is 
arbitrary this shows that \(|B_n|=0\). 

We now have that for a. a. \(x\in A\)
\[
\Big|\frac{\beta(f(x+h))-\beta(f(x))}{h}\Big| 
\le L\,\Big|\frac{f(x+h)-f(x)}{h}\Big|
\]
for \(|h|\) small enough, 
so \((\beta\circ f)'(x)=0\) and 
\(\beta'(f(x)) f'(x)=0\). This proves the 
claim of the Lemma for a. a. \(x\in(a,b)\). 
\end{proof}

\section*{Appendix B: Time regularity for the spatially 
homogeneous Boltzmann equation}

We show that the solution of the Boltzmann 
equation~\eqref{eq:bosh} under the conditions 
of Theorem~\ref{thm:1} is smooth with respect 
to time, together with its moments of any order. 

For \(k\ge 0\) 
%and positive constants \(m_0\), \(m_1\), 
%\(m_k^*\) given we introduce
%\[
%D_k = \Big\{f\in L^1 (\R^n) : 
%\ir f\,dv = m_0,
%\ir f |v|^2\, dv = m_1,
%\ir f|v|^{2k} \,dv \le m_k^*
%\Big\}. 
%\]
%We also use 
we introduce
the following weighted Lebesgue spaces
\begin{equation}
\label{eq:L1k}
L^1_k(\R^\di) = \Big\{f\in L^1 (\R^n) : 
\ir |f|\,(1+|v|^{2})^k\,dv <+\infty\Big\}
\end{equation}
with the norms defined by the integrals appearing 
in~\eqref{eq:L1k}. The regularity result that 
we used in Section~\ref{sec:mom} is the following. 

\begin{proposition}
\label{prop:treg}
Let \(f\) be the unique solution of the Boltzmann 
equation~\eqref{eq:bosh}
that preserves the total mass and energy. Assume that 
\(f_0\in L^1_k(\R^\di)\), \(k>1+\frac{\beta}{2}\). 
Then \(f\in C^1\big([0,+\infty);L^1_p(\R^\di)\big)\) 
for any \(p<k-\frac{\beta}{2}\).
\end{proposition}

The proof of Proposition~\ref{prop:treg} depends on the 
following continuity property of the nonlinear operator 
\(Q(f)\). 

\begin{lemma} 
\label{lem:Qcont}
Let the pair of positive numbers \((k,p)\)
satisfy \(k>p+\frac{\beta}{2}\). Then \(Q(f)\) is 
continuous on \(L^1_k(\R^\di)\) as a mapping 
\(L^1_k(\R^\di)\to L^1_p(\R^\di)\). Moreover, we have the following 
H\"older estimate for any \(f,g\in L^1_k(\R^\di)\)
\[
\|Q(f)-Q(g)\|_{L^1_p} \le C_{p}\Big(
\|f-g\|_{L^1}^{1-\frac{p+\frac{\beta}{2}}{k}}
+\|f-g\|_{L^1} \Big), 
\]
where the constant \(C_{p}\) depends on \(p\) and 
on the upper bound of the \(L^1_k\)-norms of \(f\) 
and \(g\).  
\end{lemma} 

\begin{proof} Using the weak form of \(Q(f)\) and \(Q(g)\) we compute
\begin{equation*}
\begin{split} 
&\ir |Q(f)-Q(g)|\, (1+|v|^{2})^p\,dv\\
& 
= \ir \ir \int_{S^{n-1}} (ff_*- gg_*)\, B (v-v_*,\sigma)\,
\Big( \text{sign}\big(Q(f)'- Q(g)'\big) (1+|v'|^{2})^p
\\
\noalign{\vskip6pt}
& 
\qquad\qquad\qquad\qquad\qquad\qquad\qquad\qquad
- \text{sign}\big(Q(f)-Q(g)\big)(1+|v|^{2})^p\Big) \, 
d\sigma\,dv\, dv_*\\
\noalign{\vskip6pt}
& \le\, 2^{p+1} \!\ir\ir 
|ff_* - gg_*|\,
|v-v_*|^\beta\, \bigl((1+|v|^{2})^p + (1+|v_*|^{2})^p\bigr)
\, dv\, dv_*\\
\end{split}
\end{equation*}
Since 
\[
\begin{split}
|v-v_*|^\beta\, (1+|v|^{2})^p & 
\le (1+|v_*|^2)^{\frac{\beta}{2}}(1+|v|^{2})^p 
+ (1+|v|^{2})^{p+\frac{\beta}{2}}
\\
& \le 2\, \bigl((1+|v|^{2})^{p+\frac{\beta}{2}} +
(1+|v_*|^{2})^{p+\frac{\beta}{2}}  \bigr)
\end{split}
\] 
and 
\(
|ff_* - gg_*| \le \frac12\,|f-g|\,|f_*+g_*| +  \frac12\,|f+g|\,|f_*- g_*|
\),
we obtain
\begin{equation*}
\begin{split}
& \|Q(f)-Q(g)\|_{L^1_p}\\
&\le 2^{p+3} \ir \ir 
|f+g|\,|f_*-g_*|\,
\bigl((1+|v|^{2})^{p+\frac{\beta}{2}} +
(1+|v_*|^{2})^{p+\frac{\beta}{2}}\bigr)\,dv\, dv_*
\\
\noalign{\vskip6pt}
&\le 2^{p+3} \,\,\|f+g\|_{L^1_k}\, \bigl(\|f-g\|_{L_{p+\frac{\beta}2}^1}
+ \|f-g\|_{L^1}\bigr).
\end{split}
\end{equation*}
We use the interpolation inequality~\eqref{eq:mint} 
with \(k_1=p+\frac{\beta}2\) to get
\begin{equation*}
\begin{split}
\|f-g\|_{L_{p+\frac{\beta}2}^1} 
& \le \|f-g\|_{L^1}^{1-\frac{p+\frac{\beta}2}k} 
\|f-g\|_{L_k^1}^{\frac{p+\frac{\beta}2}k}\\
& \le \bigl(\|f\|_{L^1_k}+\|g\|_{L^1_k}\bigr)^{\frac{p+\frac{\beta}2}k} 
\|f-g\|_{L^1}^{1-\frac{p+\frac{\beta}2}k}.
\end{split}
\end{equation*}
Substituting this bound into the previous estimate we 
obtain the H\"older estimate stated in the Lemma. This 
completes the proof. 
\end{proof}
\begin{proof}[Proof of Proposition~\ref{prop:treg}]
We fix \(T>0\). 
By the results of Arkeryd and Elmroth~\cite{Ar,El} (see part (iii) 
of Theorem~\ref{alreadyknown}),
\(f\) belongs to \(L^\infty([0,+\infty);L^1_k(\R^\di))\). By 
Lemma~\ref{lem:Qcont}, 
\begin{equation}
\label{eq:QLp}
(1+|v|^2)^p \,Q(f) \in L^1((0,T)\times\R^\di), 
\quad\text{for}\quad p<k-\frac{\beta}{2}
\end{equation}
%for \(p<k-\frac{\beta}{2}\). 
The mild form of~\eqref{eq:bosh}, together with the 
regularity condition~\eqref{eq:QLp} imply that \(f\)
is weakly differentiable and \(\partial_t f = Q(f)\) 
in the sense of distributions on \((0,T)\times\R^\di\). 
Hence,
\[
f\in W^{1,1}((0,T);L^1_p(\R^\di))
\]
and therefore (cf. \cite[p.\,286]{Ev}), 
\(
f\in C\big([0,T];L^1_p(\R^\di)\big)
\). By the continuity of \(Q(f)\) established in 
Lemma~\ref{lem:Qcont} it follows that 
\(\partial_t f\in C\big([0,T];L^1_p(\R^\di)\big)\), where 
\(\partial_t f\) is the weak time-derivative 
of \(f\). It is then easy to verify directly that 
\(f\) is strongly differentiable on \((0,T)\) with 
values in \(L^1_p(\R^\di)\), and its derivative is 
continuous on \([0,T]\). Since \(T\) is arbitrary, 
we obtain the conclusion of the Lemma. 
%Indeed, let \(g=\partial_t f\). Then 
%\[
%f(v,t+\dt)- f(v,t) = \int_t^{t+\dt} g(v,\tau)\, d\tau, 
%\] 
%for a. a. \(v\in \R^\di\). 
%We have 
%\[
%\begin{split}
%\ir \Bigl|\frac{f(v,t+\dt) - f(v,t)}{\dt}-g(v,t)\Bigr| dv
%=
%\frac{1}{\dt}\ir \Bigl|\int_t^{t+\dt} g(v,\tau)\, d\tau 
%-g(v,t)\,\dt\Bigr|\, dv
%\\
%\le \frac{1}{\dt}\ir \int_t^{t+\dt}\Bigl
%|g(\tau,v)-g(t,v)\Bigr|\, d\tau\, dv
%\\
%\ir \Bigl|\Bigr|\, dv
%\\
%\ir \Bigl|\Bigr|\, dv
%\end{split}
%\]
\end{proof}
\begin{remark}
As a consequence of Proposition~\ref{prop:treg}, if the 
moments of all orders are finite initially then they 
are continuously differentiable functions of time. 
By iterating the argument we used in the proof above 
one can show that in fact then 
\(f\in C^\infty([0,\infty);L^1_k(\R^\di))\), for 
any \(k\ge0\).   
\end{remark}

\section*{Appendix C: Carleman's representation}

\begin{lemma}
\label{lem:qp-c}
Let \(Q^+(f,g)\) be defined by \eqref{eq:bforms} 
and let \(f=f(v)\) and  \(g=g(v)\), \(v\in \R^\di\) be smooth 
functions, decaying rapidly at infinity. Then  
\begin{equation*}
Q^+(f,g)(v) = 
2^{\di-1}\ir \frac{f(v'_*)}{|v-v'_*|} \,\int_{E_{v,v'_*}} 
\frac{g(v')\,B(2v-v'-v'_*,\tfrac{v'-v'_*}{|v'-v'_*|})}{|v'-v'_*|^{\di-2}}\,
d\pi_{v'}\,dv'_*, 
\end{equation*}
where 
\(E_{v,v'_*}\) is the hyperplane
\(\{v'\in \R^\di\;|\;(v'-v)\cdot(v'_*-v) = 0\}\)
and \(d\pi_{v'}\) denotes the Lebesgue measure on this 
hyperplane.
\end{lemma}
\begin{proof} 
Using the change of variables \(u=v-v_*\), and recalling 
the definition of the delta function of a quadratic 
form, see~\cite{GeSh1}, we have 
\begin{equation}
\label{eq:qpdelta}
Q^+(f,g)(v) = \ir\ir f(v'_*) \,g(v')\,B(u,k)\,
\delta\bigl(\tfrac{|k|^2-1}{2}\bigr)\,dk\,du,
\end{equation}
where \(v' = v - u + \frac{1}{2}\,(u+|u| k )\) and
\(v'_* = v - \frac{1}{2}\,(u+|u|k )\). 
We further set 
\(
z=-\frac{1}{2}(u+|u|k)
\); 
for every \(u\) fixed this defines a linear map \(k\mapsto z\) 
with determinant \(\bigl(\frac{|u|}{2}\bigr)^\di\). 
We also have  
\[
k = -\frac{2z+u}{|u|}
\quad\text{and}\quad
\frac{|k|^2-1}{2} 
=\frac{|2z+u|^2-|u|^2}{2|u|^2}
=\frac{2z\cdot(z+u)}{|u|^2}.  
\]
With this change of variables the integral 
in~\eqref{eq:qpdelta} can be written as
\[
\ir\ir 
\bigl(\tfrac{2}{|u|}\bigr)^{\di}
f(v+z)\,g(v-u-z)\,B(u,-\tfrac{2z+u}{|u|})\,
\delta\bigl(\tfrac{2z\cdot(z+u)}{|u|^2}\bigr)\,dz\,du.
\]
We set \(y=-z-u\); then \(|u|=|y+z|\) and 
\(\delta\bigl(\tfrac{2z\cdot (z+u)}{|u|^2}\bigr)
=\tfrac{|y+z|^2}{2}\,\delta (z\cdot y )
\). Further, for any test function \(\ph\),
\[ 
\ir\delta(z\cdot y)\,\ph(y)\,dy 
= |z|^{-1} \int_{z\cdot y=0} \ph(y)\,d\pi_y,        
\]
where \(d\pi_y\) is the Lebesgue measure on 
the hyperplane \(\{y:z\cdot y=0\}\). This 
yields
\[
\begin{split}
& Q^+(f,g) (v)
\\
& = {2^{\di-1}}\!\int_{z\in\R^\di} 
\int_{y\cdot z=0} f(v+z)g(v+y)\, |z|^{-1}\,
\!|y+z|^{n-2}\,B(-y-z,\tfrac{y-z}{|y+z|})\,d\pi_y\, dz
\end{split}
\]
We now return to the original notations
\(v'_* = v+z\), \(v' = v+y\) and perform 
the corresponding changes of variables to obtain 
the expression for \(Q^+(f,g)\) stated in the Lemma. 
\end{proof}

\begin{remark}
The above result takes a particularly simple form 
in the case of the hard-sphere model in \(\R^3\); 
in that case \(B(v-v_*,\sigma)=\frac{1}{4\pi}|v-v_*|\) 
and 
\[
Q^+(f,g) (v)
= \int_{\R^3} \frac{f(v'_*)}{\pi |v-v'|}
\int_{E_{v,v'_*}} \!g(v')\,
d\pi_{v'} dv'_*. 
\]
\end{remark}
\bigskip

\Ack{The research of the first author was partially supported 
by NSF under grant DMS-0507038. The  second author was partially 
supported by PIMS and by NSERC under operating grant 7847. 
The third author acknowledges support from the HYKE European 
network, contract HPRN-CT-2002-00282. Support from the Institute 
from Computational Engineering and Sciences at the University 
of Texas at Austin is also gratefully acknowledged.}

\bibliographystyle{acm}
%\bibliography{umb}

\signig
\signvp
\signcv

\end{document}